\pdfoutput=1
\documentclass[journal]{IEEEtran}
\usepackage{cite}

\ifCLASSINFOpdf
  \usepackage[pdftex]{graphicx}
\else
\fi
\usepackage[cmex10]{amsmath}
\usepackage{amssymb}  

\usepackage{algorithm}
\usepackage{algpseudocode}

\ifCLASSOPTIONcompsoc
 \usepackage[caption=false,font=normalsize,labelfont=sf,textfont=sf]{subfig}
\else
 \usepackage[caption=false,font=footnotesize]{subfig}
\fi
\newtheorem{thm}{Theorem}
\newtheorem{lem}[thm]{Lemma}

\newtheorem{prop}[thm]{Proposition}

\newtheorem{assum}{Assumption}

\newtheorem{rem}{Remark}

\usepackage{soul,color}
\usepackage{mydefs}
\def\alert #1{\textcolor{black}{#1}}

\makeatletter
\newcommand{\removelatexerror}{\let\@latex@error\@gobble}
\makeatother

\begin{document}
%
\title{A Bregman Splitting Algorithm for Distributed Optimization over Networks$^*$ \thanks{$^*$Submitted initially to the editors on October 20, 2015 and currently revised on August 28, 2016. Some partial result has been submitted to the 55th IEEE Conference on Decision and Control for possible presentations.}}
%
%
%

\author{Jinming~Xu,~\IEEEmembership{Student Member,~IEEE,}
        Shanying~Zhu,~\IEEEmembership{Member,~IEEE,}
        Yeng Chai~Soh,~\IEEEmembership{Senior Member,~IEEE,}
        and~Lihua~Xie,~\IEEEmembership{Fellow,~IEEE}}
\maketitle

\begin{abstract}
We consider distributed optimization problems in which a group of agents are to collaboratively seek the global optimum through peer-to-peer communication networks. The problem arises in various application areas, such as resource allocation, sensor fusion and distributed learning. We propose a general efficient distributed algorithm--termed Distributed Forward-Backward Bregman Splitting (D-FBBS)--to simultaneously solve the above primal problem as well as its dual based on Bregman method and operator splitting. The proposed algorithm allows agents to communicate asynchronously and thus lends itself to stochastic networks. This algorithm belongs to the family of general proximal point algorithms and is shown to have close connections with some existing well-known algorithms when dealing with fixed networks. However, we will show that it is generally different from the existing ones due to its effectiveness in handling stochastic networks. With proper assumptions, we establish a non-ergodic convergence rate of $o(1/k)$ in terms of fixed point residuals over fixed networks both for D-FBBS and its inexact version (ID-FBBS) which is more computationally efficient and an ergodic convergence rate of $O(1/k)$ for D-FBBS over stochastic networks respectively. We also apply the proposed algorithm to sensor fusion problems to show its superior performance compared to existing methods.
\end{abstract}

%
\IEEEpeerreviewmaketitle

\section{Introduction}
%
%
%
%

 \IEEEPARstart{D}{istributed} optimization has received much attention recently from the control and machine learning communities due to its wide applications. Examples include formation control~\cite{DynamicOptim_formation_flight}, resource allocation~\cite{Dual_decomposition_M.Johansson_S.Boyd}, information fusion in wireless sensor networks~\cite{Rabbat2004_WSN,chen2012diffusion,zhu2013distributed} and distributed learning~\cite{Predd2009_DistributedLearning,mateos2010distributed} and just to name a few. Common features to these examples are that the system typically has a large number of agents (e.g., sensors) involved without any centralized coordinator and that resources, such as sensing, communication and computation, are usually scattered throughout the network. These features naturally lead to the necessity of completely distributed algorithms that can operate merely based on local information and are robust against the changes of network topology. 

 The existing literatures on distributed optimization can be generally categorized into two main streams: consensus-based approach and dual-decomposition-based approach. In the context of distributed computation, consensus theory is particularly suitable for distributed implementation of algorithms as it allows agents to obtain global results by merely taking local actions even over time-varying networks~\cite{Consensus_Olfati-Saber}. Specifically, the study of consensus-based approach can be dated back to the seminal work~\cite{tsitsiklis1986distributed}, where Tsitsiklis and Bertsekas et al. analyzed distributed gradient-like optimization algorithms in which a number of processors perform computations and exchange messages asynchronously intending to minimize a certain cost function.  In line with this work, Nedic and Ozdaglar applied consensus theory into more general optimization problems where each agent has its own local private cost function and aims to minimize the sum of cost functions through cooperation with neighbors over networks which can be time-varying~\cite{Subgradient_Consensus_Nedic}. To avoid the worst-case bounded communication assumption as usually made by the abovementioned works, Lobel and Ozdaglar~\cite{lobel2011distributed} considered the same algorithm for stochastic networks where communication links are subject to random failures. In~\cite{Duchi2012_DualAveraging}, an interesting dual subgradient averaging method using proximal function is proposed for solving the same problem which is shown to have better convergence results in terms of network scaling and is claimed to be able to deal with stochastic networks. Some extensions have also been made for constrained cases~\cite{Nedic2010_ConstrainedConsensus,Zhu2012_DistrOptiConstraint}, situations where only noisy gradient is available~\cite{Srivastava2011_StochasticConstraintOpti,Nedic2014_LipschitzGrad}, and directed graphs~\cite{Nedic2015_Pushsum}. Most algorithms as abovementioned require decaying stepsizes and bounded (sub)gradient assumption for achieving the exact optimum while using constant stepsize can only guarantee its convergence to a neighborhood of the same order of the stepsize~\cite{yuan2013convergence}. However, using diminishing stepsize always results in weaker convergence rates. In particular, both the dual-averaging method~\cite{Duchi2012_DualAveraging} and (sub)gradient-push method~\cite{Nedic2015_Pushsum} show a convergence rate of $O(\ln k/\sqrt{k})$ in terms of \emph{objective error} with a decaying stepsize of $1/\sqrt{k}$ in the ergodic sense.  The convergence rate is improved to $O(\ln k/k)$ in~\cite{Jakovetic2014_FastGradMethod} by employing the accelerated \emph{Nesterov} method and also in~\cite{Nedic2014_LipschitzGrad} for strongly convex functions with Lipschitz gradients.

On the other hand, dual decomposition also lends itself to distributed implementation of algorithms for large-scale optimization problems which are separable in the dual domain~\cite{Dual_decomposition_R.Murray,Book:ConvexOptim}. It has been shown that distributed optimization problems can be transferred to an equivalent constrained optimization problem~\cite{Cortes2014_conststepsize} for which we have a lot of existing techniques available. In particular, doing this allows us to transfer the problem as a saddle point problem for which the traditional Arrow-Hurwiz-Uzawa method can be applied~\cite{benzi2005numerical}, especially the augmented Lagrangian method which permits better convergence performance~\cite{fortin2000augmented}. However, introducing augmented term results in coupling issues among cost functions. To overcome this, a popular approach called as alternating direction method of multiplier (ADMM) is proposed which is shown to have very good convergence performance even for large-scale problems~\cite{Boyd2011_DistOpti}. \alert{Building on this method, distributed versions of ADMM~\cite{Mota2013_D-ADMM,wei2012distributed_admm,shi2014linear}, along with its inexact versions~\cite{chang2015multi,ling2015dlm}, are developed for solving the distributed optimization problem by explicitly dealing with the coupling topology. Following the same idea, some attempts are being recently made to improve the convergence performance \cite{mokhtari2016decentralized} and for generalization \cite{hong2015stochastic}.} It is well known that the abovementioned methods are specific applications of proximal point algorithm and operator splitting~\cite{rockafellar1976monotone,eckstein1992douglas,bauschke2011convex}, which have been widely applied in signal processing~\cite{combettes2011proximal} and image processing\cite{chambolle2011first}.  Our algorithm also belongs to this big family. The advantage of using dual-decomposition-based approach is that it is able to achieve the exact optimum and generally obtain an ergodic convergence rate of $O(1/k)$. However, differ from the consensus-based approach, the effectiveness of the dual-decomposition-based approach relies heavily on the knowledge of the (coupling) structure of the problem. Thus, most algorithms\footnote{Randomized iteration of algorithms requires the weight matrix to be constant~\cite{jakovetic2015linear,bianchi2014primal,hong2015stochastic}. Thus it is deemed applicable only to fixed networks.} developed under this framework are only capable of dealing with fixed networks where the coupling topology is not varying with time and its application to stochastic networks is still unclear due to the difficulty in analyzing the convergence properties. 
\par
 \alert{In this paper, we propose a general efficient distributed algorithm based on Bregman method and operator splitting to solve distributed optimization problems over both fixed and stochastic networks. The algorithm, termed Distributed Forward-Backward Bregman Splitting (D-FBBS), is able to solve both the primal and dual problems simultaneously (cf.~Section~\ref{sec:prob_description}). We also show that the proposed algorithm is general and has close connections with some existing well-known algorithms. Indeed, the Bregman splitting method used to develop this algorithm provides a unified framework which allows us to easily recover most existing algorithms such as P-EXTRA/EXTRA~\cite{shi2015extra,shi2015proximal}, D-ADMM/DLM~\cite{shi2014linear,ling2015dlm} and C-ADMM/IC-ADMM~\cite{mateos2010distributed,chang2015multi} by properly designing the consensus constraint\footnote{Similar with ADMM-based framework, one can further exploit this structure to broaden the class of distributed algorithms for particular problems.} and employing certain splitting schemes (cf. Section \ref{sec:variations_connections}). To the best of our knowledge, we are the first to introduce the Bregman splitting method to solve distributed optimization problems and, most importantly, show its effectiveness in dealing with stochastic networks (i.e., using Bregman distance as a common Lyapunov function). In particular, with proper assumptions, we establish a non-ergodic convergence rate of $o(1/k)$ and an ergodic convergence rate of $O(1/k)$ in terms of fixed point residuals (FPR) for fixed networks and stochastic networks respectively,  which are the best known convergence rates in the existing literature. For more computational efficiency, we also provide an inexact version of D-FBBS (termed ID-FBBS) for cost functions that have Lipschitz gradients and, similarly, show its convergence with a non-ergodic rate of $o(1/k)$ for fixed networks. It should also be noted that our formulation is node-based and is thus naturally suitable for distributed computation unlike the edge-based counterparts~\cite{wei2012distributed_admm,bianchi2014primal,hong2015stochastic}.} 

Our work is closely related to some existing algorithms. In particular, Jakovetic at al.~\cite{jakovetic2015linear} studied the linear convergence rate of a class of distributed augmented Lagrangian (AL) algorithms for twice continuously differentiable cost functions with a bounded Hessian when there are sufficient inner iterations of consensus being carried out. In the recent work~\cite{shi2015extra}, Shi et al. proposed a similar algorithm termed EXTRA, where the cost function is assumed to have Lipschitz gradients. When the cost function is also strongly convex, linear convergent rate can be achieved as well. For generalization, they extend the algorithm to general convex problems and composite convex problems, yielding P-EXTRA and PG-EXTRA respectively~\cite{shi2015proximal}. Our algorithm, though is generally different in its nature, has close connections with them in the sense that, with proper parameter settings,  these existing algorithms can be shown to be equivalent to ours with certain splitting schemes (cf.~Section~\ref{sec:variations_connections}). \alert{Note that the Bregman splitting method has been used in developing the Bregman Operator Splitting (BOS) algorithm~\cite{zhang2010bregmanized} which deals with general equality-constrained problems.  However, our algorithm differs from BOS in that we deal with distributed optimization problems (i.e., focusing on distributed implementation) and consider an asymmetric saddle point problem (cf. Eq.~(\ref{eq:asymmetric_saddlepoint_probs})), which, as we will see later, allows us to effectively deal with stochastic networks.}

{\it Notation:}  A vector is regarded as a column vector unless clearly stated otherwise. We denote by $x_i$ the $i$-th component of a vector $x$ and $x_{i,k}$ the $i$-th component of a vector $x$ at time $k$. 
We use $\bar{x}=\avector x$ to denote the vector whose elements are the average of all components of $x$. In addition, we use $\otimes$ to denote the Kronecker product of two matrices, `$\ones$' all-ones column vector with proper dimension, $I$ the identity operator or identity matrix with proper dimension, $\mathcal{H}$ the Euclidean space, $\innprod{\cdot}{\cdot}$ the inner product and $\Gamma_0(\mathcal{H})$ the class of proper lower semi-continuous convex functions from $\mathcal{H}$ to $(-\infty,+\infty]$. Moreover, we use  $O(\cdot)$ to denote a sequence having an ergodic convergence rate if the rate is stated in terms of the running average $\hat{x}_k=\frac{1}{k}\sum^{k-1}_{i=0}x_k$ and $o(\cdot)$ the non-ergodic convergence rate if stated in terms of $x_k$.

\section{Preliminaries}
\subsection{Subdifferential and Monotone Operator}
An operator $T$ on an Euclidean Space $\mathcal{H}$ is a set-valued mapping, i.e., $T:\mathcal{H}\rightarrow 2^\mathcal{H}$. \alert{The subdifferential of $f:\mathcal{H}\rightarrow (-\infty,+\infty]$ at a point $x$ is the set-valued operator $$\partial f:=\{p\in\mathcal{H}|f(y)\geq f(x)+\innprod{p}{y-x},\forall y\in\mathcal{H}\}.$$}
The graph of $T$ is defined as $\text{gra}~T=\{(x,y)\in\mathcal{H}\times\mathcal{H}|y\in Tx\}$. The operator $T$ is called monotone if $$\forall (x,y), (x',y')\in \text{gra}~T ~~~~ \innprod{x-x'}{y-y'}\geq 0.$$
A monotone operator is said to be maximal if there is no monotone operator $T'$ such that $\text{gra}~T\subset\text{gra}~T'$, e.g., the subdifferential $\partial f$ of $f\in\Gamma_0(\mathcal{H})$ is maximally monotone.
The resolvent of $T$ is defined as ${ R}_{\tau T}=(I+\tau T)^{-1}$. The proximity operator of a convex functional $f$ is just an application of resolvent operator with $\tau> 0$ and is defined as follows $${\bf prox}_{\tau f}(v)=\underset{x\in\mathcal{H}}{\text{arg}\min}~\{f(x)+\frac{1}{2\tau}\norm{x-v}^2\}.$$


\subsection{Bregman Distance and $G$-space}
Consider a convex functional $f:\mathcal{H}\rightarrow\mathcal{R}\cup \{\pm\infty\}$, the generalized Bregman distance between $x\in\mathcal{H}$ and $x'\in\mathcal{H}$ is defined as (cf.~see \cite{osher2005iterative} for more details)
$$D^q_f(x,x')=f(x)-f(x')-\innprod{q}{x-x'},$$
where $q\in\partial f(x')$ is a subgradient evaluated at $x'$. Bregman distance has the following properties:
\begin{itemize}
\item $D^q_f(x,y)\geq 0$,
\item $D^q_f(y,x)\geq D^q_f(z,x)$ if $z\in[x,y]$.
\end{itemize}

Given a symmetric positive definite matrix $G$, we define a $G$-space as well as its induced norm as follows:
$$\innprod{x}{x'}_G=\innprod{Gx}{x'}~\text{and}~\norm{x}_G=\sqrt{\innprod{Gx}{x}},~\forall x,x'\in\mathcal{H}.$$
\subsection{Saddle Point and Fenchel Duality}
A pair $(x^\star,y^\star)$ is said to constitute a saddle point of $\psi$ in the domain $\mathcal{D}$ if the following condition holds 
\begin{equation}\label{eq:saddlepoint_definition}
\psi(x^\star,y)\leq\psi(x^\star,y^\star)\leq\psi(x,y^\star), \forall (x,y)\in \mathcal{D},
\end{equation}
where $\psi(x,y)$ refers to the Lagrangian in the sequel.

Let $f:\mathcal{H}\rightarrow\mathcal{R}\cup \{\pm\infty\}$ be a convex functional. Then, its convex conjugate (Fenchel's dual) $f^\ast$ is defined as follows: 
$$f^\ast(y):=\sup_{x\in \mathcal{H}}\Bracket{\innprod{x}{y}-f(x)},$$
where $y$ is the dual variable. Note that $f^{**}=f$ by Fenchel–-Moreau theorem if $f\in\Gamma_0(\mathcal{H})$.
\section{Problem Statement and Assumptions}
In this section, we state our problem and the network model we consider in this paper. In addition, we formulate the distributed optimization problem as a \emph{primal-dual} (saddle-point) problem and show that by solving the \emph{primal-dual} problem (i.e, attaining the saddle point), we will not only obtain the optimal solution to the primal distributed optimization problem but also the optimal solution to its dual problem which is also an interesting problem to be solved separately.
\subsection{Network Model}
We consider a communication network represented by an undirected graph $\mathcal{G}=(\mathcal{V},\mathcal{E})$,  where each vertex $v_i\in\mathcal{V}$ denotes an agent and each edge $e_{ij}=(v_i,v_j)\in\mathcal{E}$ the communication link which is assigned with a \emph{positive} weight~$w_{ij}$. We use $\mathcal{N}_i:=\{j|e_{ij}\in \mathcal{E}\}$ to denote all the neighbors of agent $i$. For stochastic scenarios, we model each communication link as a probabilistic event. In particular, we assume each link fails according to certain i.i.d. random process. As a result, the weight matrix $W=\{w_{ij}\}$ associated with the network can be sampled from a probability space~$\mathcal{F}=(\Omega,\mathcal{B},\mathcal{P})$, where $\Omega$ is the set of all feasible stochastic matrix $W\in \mathcal{R}^{m\times m}$, $\mathcal{B}$ the Borel $\sigma$-algebra on $\Omega$ and $\mathcal{P}$ the probability measure. It should be noted that this network  model includes as a special case the gossip algorithm~\cite{boyd2006randomized} and allows the interaction topology among agents changing with time.


For clarity of presentation, let us first consider a fixed network for which we make the following assumption:
\begin{assum}\label{assum:fixednet_weightmatrix}
The \emph{non-negative} weight matrix $W$ as defined above satisfies the following conditions:
\begin{itemize}
\item Positive-definite: $W^T=W$ and $W>0$,
\item Stochasticity: $W\ones=\ones$ or $\ones^TW=\ones^T$,
\item Connectivity: $\rho\left(W-\avector\right)<1$.
\end{itemize}
where $\rho(\cdot)$ is the spectral radius of a given matrix.
\end{assum}
\begin{rem}\label{rem:weightmatrix_eigen_speed}
\alert{Based on the above assumption, we immediately know that $W$ has a simple eigenvalue one and all the other eigenvalues $\lambda_i(W)\in (0,1)$, and $null(I_m-W)=span\{\mathbf{1}\}$.  Also, note that the weight matrix has significant impacts on the speed of information diffusion and thus needs to be properly designed for better convergence performance~\cite{xiao2004fast,kar2008sensor}.}
\end{rem}
\subsection{Distributed Optimization Problem (DOP)}\label{sec:prob_description}
Consider a network consisting of $m$ agents as described above. The objective of the network is to minimize the following function in a cooperative way: 
\begin{equation}
F(\theta)=\sum_{i=1}^mf_i(\theta), \tag{DOP}
\end{equation} where $\theta\in\mathcal{R}^d$ is the global variable to be optimized while $f_i:\mathcal{R}^d\rightarrow\mathcal{R}\cup \{\pm\infty\}$ is the local cost function available only to agent $i$, for which we make the following assumption:  
\begin{assum}\label{assum:costfuctions_proper_closed_convex}
The cost functions are proper, closed and convex, i.e., $f_i\in\Gamma_0(\mathcal{H}), \forall i\in \mathcal{V}$. 
\end{assum}

Let $f(x):=\sum_{i=1}^mf_i(x_i)$, where $x_i\in\mathcal{R}^d$ is the locally stored variable of agent~$i$ for the optimal solution of the DOP problem and $x=[x_1^T,x_2^T,...,x_m^T]^T\in\mathcal{R}^{md}$ the collection of all local variables. Then, it follows from Assumption~\ref{assum:costfuctions_proper_closed_convex} that $f$ is proper, closed and convex, i.e., $f\in\Gamma_0(\mathcal{H})$. In addition, defining a consensus space $\cspace:=\{[\theta^T,\theta^T,...,\theta^T]^T|\theta\in\mathcal{R}^d\}$, then the DOP problem can be rewritten as follows 
\begin{equation}\label{prob:dops_indicatorfun}
\min_{x\in\mathcal{R}^{md}} f(x)+\iota_\cspace(x),
\end{equation}
where $\iota_\cspace(x)$ is the indicator function defined as follows
\begin{equation*}\label{eq:indicatorfun}
\iota_\cspace(x)=
\begin{cases}
0,~&\text{if}~x\in \cspace,\\
\infty,~&\text{otherwise.}
\end{cases}
\end{equation*}
Similar with~\cite[Lem.~3.1]{Cortes2014_conststepsize}, since $null\{I_m-W\}=span\{\ones\}$ by Assum. \ref{assum:fixednet_weightmatrix} (cf.~Rem.~\ref{rem:weightmatrix_eigen_speed}), the problem (\ref{prob:dops_indicatorfun})  is equivalent to the following problem, termed optimal consensus problem~\ocp:
\begin{equation}\label{prob:dops}
\min_{x\in\mathcal{R}^{md}} f(x)=\sum_{i=1}^mf_i(x_i)~~~{\bf s.t.}~(I_{md}-W\otimes I_d)x=0, 
\end{equation}
 where $I_{md}=I_m\otimes I_d$ is the identity matrix\footnote{Throughout the paper, we will, once it is clear, suppress the subscript of certain variables for brevity.}. 

Correspondingly, the dual formulation of the problem (\ref{prob:dops_indicatorfun}) is as follows~\cite[Def. 15.10]{bauschke2011convex}:
\begin{equation}\label{prob:dops_dual_indicatorfun}
\min_{y\in\mathcal{R}^{md}} f^\ast(y)+\iota^\ast_\cspace(-y),
\end{equation}
where $f^\ast$ and $\iota^\ast_\cspace(\cdot)$ are the convex conjugates of $f$ and $\iota_\cspace(\cdot)$ respectively. Likewise, since $\iota^\ast_\cspace(\cdot)$ indicates the orthogonal space of $\cspace$ (denoted as $\cspace^\perp$) the problem (\ref{prob:dops_dual_indicatorfun}) is equivalent to the following optimal exchange problem~\oep:
\begin{equation}\label{prob:dops_dual}
\min_{y\in\mathcal{R}^{md}} f^\ast(y)=\sum_{i=1}^mf^\ast_i(y_i)~~~{\bf s.t.}~(\ones\otimes I_d)^Ty=0,
\end{equation}
where $y=[y_1,y_2,...,y_m]^T\in\mathcal{R}^{md}$ is the dual variable. 

The above analysis shows that the OCP problem (primal) and the OEP problem (dual) are, in fact, dual to each other. As we will show later, the following proposed distributed algorithm can solve both the primal and dual problem simultaneously. In other words, it provides an alternative \emph{distributed} way to solve the OEP problem which usually needs to be solved by centralized or parallel approaches~\cite{Boyd2011_DistOpti}. Note that, in what follows, the above OCP and OEP problems will be termed together as a \emph{primal-dual} problem.  

We introduce the following augmented Lagrangian associated with the above \emph{primal-dual} problem.
\begin{equation}\label{eq:dops_lagrangian}
\begin{aligned}
\psi(x,y)=f(x)-y^Tx+\frac{1}{2\gamma}\norm{x}^2_{I_{md}-W\otimes I_d}.\\
\end{aligned}
\end{equation}
Note that the dual variable (Lagrange multiplier) $y$ plays a key role in reconciling the discrepancy of the interests of different agents for achieving global optimum. 
Without loss of generality, we assume all agents use the same parameter  $\gamma$. 

Accordingly, the KKT conditions for optimality of the above \emph{primal-dual} problem can be depicted as follows:
\begin{subequations}\label{eq:dops_optimality_conditions}
\begin{align}
&(\text{Primal Feasibility})~~~~~~~~~(I_{md}-W\otimes I_d)x^\star=0, \label{eq:dops_optimality_conditions_primalfeasible}\\
 &(\text{Dual Feasibility})~~~~~~~~~~~~~~~~(\ones\otimes I_d)^Ty^\star=0, \label{eq:dops_optimality_conditions_dualfeasible}\\
 &(\text{Lagrangian Optimality})~~~~~~~~~~~y^\star\in\partial f(x^\star). \label{eq:dops_optimality_conditions_Lagrangeoptimal}
 \end{align}
\end{subequations}
As we will show later, the dual feasibility~(\ref{eq:dops_optimality_conditions_dualfeasible}) can be always guaranteed (cf.~Lem.~\ref{lem:ConservProp}) by the following proposed algorithm with proper initialization. Thus, in the sequel, we will restrict our attention to the domain $\mathcal{D}=\{(x,y)|x\in\mathcal{R}^{md},~y \in \cspace^\perp\}$.

For the problem to be feasible and the strong duality holds, we make the following assumption:

\begin{assum}\label{assum:existence_saddlepoint}
The augmented Lagrangian $\psi$, as defined in~(\ref{eq:dops_lagrangian}), has a saddle point $(x^\star,y^\star)\in\mathcal{D}$.
\end{assum}

The following lemma shows the equivalence between the saddle point of~(\ref{eq:dops_lagrangian}) and the optimality conditions~(\ref{eq:dops_optimality_conditions}).
\begin{prop}[Optimality as Saddle Points]\label{prop:optim_iff_saddlepoint}
Consider a pair $(x^\star,y^\star)\in\mathcal{D}$. Then, it is a saddle point of the Lagrangian~(\ref{eq:dops_lagrangian}) if and only if it satisfies the optimality conditions~(\ref{eq:dops_optimality_conditions}). Moreover, if it is a saddle point, then $x^\star$ solves the OCP problem and $y^\star$ solves the OEP problem respectively. 
\end{prop}
\begin{IEEEproof}
See Appendix~\ref{sec:appendix}.
\end{IEEEproof}
\begin{rem}
In the rest of the paper, for simplicity and brevity, we will assume $d=1$ and the results developed in this paper can be extended to the case $d\neq 1$ without much effort. As a result, we have $\ones\otimes I_d=\ones, W\otimes I_d=W,~\mathcal{C}=span\{\ones\}$.
\end{rem}




\section{A Distributed Forward-Backward Bregman Splitting Algorithm}

In this section, we present a distributed algorithm--Distributed Forward-Backward Bregman Splitting (D-FBBS)--to solve the above \emph{primal-dual} problem. In particular, we first review the well known distributed subgradient method (DSM) and its limitation. Then, the Bregman Iterative Regularization is introduced as the basis of developing our algorithm. Further, we use the forward-backward splitting technique to split the optimization problem, leading to a separated one which can be solved efficiently in a distributed way with much cheaper computation and less communication over the network. 
\subsection{Review of Some Basic Techniques}
\subsubsection{Distributed Subgradient Methods (DSM)}
In~\cite{Subgradient_Consensus_Nedic}, Nedic and Ozdaglar proposed the following well-known distributed algorithm\footnote{To facilitate the analysis of its limitation, we only consider the case where the weight matrix is fixed and constant stepsize is employed.} to solve the OCP problem:
\begin{equation}\label{alg:DSM}
x_{k+1}=Wx_k-\gamma\cdot s_k,
\end{equation}
where $s_k\in \partial f(x_k)$ is a subgradient evaluated at $x_k$. It is easy to show that the limit point of the above algorithm, if exists, solves the inclusion: $0\in (I-W)x+\gamma \partial f(x)$, corresponding to the following minimization problem:
$$\min_xf(x)+\frac{1}{2\gamma}\norm{x}^2_{I-W},$$
which resembles penalty-based approaches that approximately solve the OCP problem. Thus, only when $\gamma$ is sufficiently small will $x$ converge to the neighborhood of the true solution. In order to obtain the exact optimum, one may use decaying stepsizes but doing so will lead to a slower convergence rate, e.g., a decaying stepsize of $\gamma_k=\frac{1}{\sqrt{k}}$ resulting in an ergodic rate of $O(\frac{\ln k}{k})$ in terms of objective error~\cite{Duchi2012_DualAveraging}. 
\subsubsection{Bregman Iterative Regularization}
In~\cite{osher2005iterative}, a Bregman-based method is introduced and revealed to be very efficient in image processing. Bregman iterative regularization method attempts to solve the following optimization problem $$\min_x J(x)+H(x)$$ iteratively by the following algorithm
\begin{equation}\label{alg:Bregman_iter_regular}
\begin{aligned}
x_{k+1}&\in\text{arg}\underset{x}{\min}~\bracket{D^{y_k}_J(x,x_k)+H(x)}\\
y_{k+1}&=y_k-\nabla H(x_{k+1}), \\
\end{aligned}
\end{equation}
where $J$ is a convex regularization functional, $H$ is a convex differentiable fitting functional in inverse problems, e.g., $H(x)=\norm{Ax-b}^2$ for linear problems, and $D^{y_k}_J(x,x_k)$ is the Bregman distance between $x$ and $x_k$ while $y_k\in\partial J(x_k)$ is the corresponding subgradient of $J$ at $x_k$. It was shown in~\cite{yin2008bregman} that the above algorithm is, in fact, equivalent to the augmented Lagrangian method for linear problems and solves the following equality-constrained problem:
\begin{equation}\label{prob:equality_const_minimum}
\min_{x} J(x) ~~~~~{\bf s.t.}~\nabla H(x)=0.
\end{equation}
Inspired from the above analysis, following the algorithm~(\ref{alg:Bregman_iter_regular}), we can easily come up with an algorithm for solving the OCP problem over fixed networks as follows:
\begin{subequations}\label{alg:Bregman_iter_regular_dops_fixednet}
\begin{align}
x_{k+1}&\in\underset{x\in\mathcal{R}^{m}}{\text{arg}\min}~\bracket{D^{y_k}_f(x,x_k)+\frac{1}{2\gamma}\norm{x}^2_{I-W}}, \label{alg:Bregman_iter_regular_dops_fixednet_x}\\ 
y_{k+1}&=y_k-\frac{1}{\gamma}(I-W) x_{k+1}. \label{alg:Bregman_iter_regular_dops_fixednet_y}
\end{align}
\end{subequations}
\begin{rem}\label{rem:Bregman_iter_regular_issues}
The above algorithm cannot be carried out in a distributed way since the $x$-update requires the linear operator $W$ to be evaluated implicitly thus requiring either the global knowledge of the network or infinite inner-loops of consensus to solve the subproblem~(\ref{alg:Bregman_iter_regular_dops_fixednet_x}), which is not practical.
\end{rem} 
\subsubsection{Forward-Backward Splitting (FBS)}
In composite optimization problems, one always needs to solve the following inclusion: $0\in (A+B)z$, where $A$ and $B$ are two operators. It is usually expensive to solve two operators together and a more efficient way is to split them into separated parts each of which is relatively easier to be evaluated. There are several splitting techniques proposed in the existing literature. Forward-Backward splitting is the one dedicated for inclusion problems where $A$ is maximally monotone and $B$ is cocoercive~\cite{bauschke2011convex}. In particular, the above inclusion can be rewritten in a split form as: $0\in (I+\tau A)z-(I-\tau B)z$, which leads to the following forward-backward splitting algorithm:
\begin{equation}\label{eq:forbacksplit_iterseq}
z_{k+1}={R}_{\tau A}(I-\tau B)z_k,
\end{equation} 
where ${R}_{\tau A}$ is the resolvent of $A$ which amounts to ${\bf prox}_{\tau f}$ when $A=\partial f$ with $f\in\Gamma_0(\mathcal{H})$. 
\subsection{Proposed Distributed Algorithm}
As abovementioned, the algorithm~(\ref{alg:Bregman_iter_regular_dops_fixednet}) does not permit for distributed implementation (cf. Remark~\ref{rem:Bregman_iter_regular_issues}). To see this, let us consider the $x$-update step~(\ref{alg:Bregman_iter_regular_dops_fixednet_x}). By the necessary condition of optimality we have 
\begin{equation}\label{eq:Bregman_iter_regular_dops_fixednet_x_equi_inclusion}
\begin{aligned}
\gamma y_k\in(I-W+\gamma\partial f) (x_{k+1}).
\end{aligned}
\end{equation} 
It is clear that the above inclusion involves computing the inverse of $W$, suffering the aforementioned issues (cf. Rem.~\ref{rem:Bregman_iter_regular_issues}). In order to solve it in a distributed way, we propose a forward-backward splitting approach as follows (cf.~Eq.~(\ref{eq:forbacksplit_iterseq})):
\begin{equation}\label{eq:Bregman_iter_regular_dops_fixednet_x_forbacksplit_inclusion}
\begin{aligned}
x_{k+1}={\bf prox}_{\gamma f} (Wx_k+\gamma y_k).
\end{aligned}
\end{equation} 
According to forward-backward splitting\cite{bauschke2011convex}, given certain $y_k$, (\ref{eq:Bregman_iter_regular_dops_fixednet_x_forbacksplit_inclusion}) is equivalent to (\ref{eq:Bregman_iter_regular_dops_fixednet_x_equi_inclusion}) only when it runs infinite steps. However, we will show that the OCP problem can still be solved when (\ref{eq:Bregman_iter_regular_dops_fixednet_x_forbacksplit_inclusion}) is executed only once per each iteration, which yields the following algorithm:
\begin{subequations}\label{alg:Bregman_iter_regular_dops_fixednet_forbacksplit_equi_seq}
\begin{align}
\gamma y_k-(x_{k+1}-Wx_k)&\in\gamma\partial f(x_{k+1}) \label{alg:Bregman_iter_regular_dops_fixednet_forbacksplit_equi_seq_x}\\
(I-W) x_{k+1}+\gamma(y_{k+1}-y_k)&= 0, \label{alg:Bregman_iter_regular_dops_fixednet_forbacksplit_equi_seq_y}
\end{align}
\end{subequations} 
or, equivalently,
\begin{subequations}\label{alg:Bregman_iter_regular_dops_fixednet_forbacksplit}
\begin{align}
x_{k+1}&=\underset{x\in\mathcal{R}^{m}}{\text{arg}\min}~\bracket{D^{y_k}_f(x,x_k)+\frac{1}{2\gamma}\norm{x-Wx_k}^2}, \label{alg:Bregman_iter_regular_dops_fixednet_forbacksplit_x}\\
y_{k+1}&=y_k-\frac{1}{\gamma}(I-W) x_{k+1}. \label{alg:Bregman_iter_regular_dops_fixednet_forbacksplit_y}
\end{align}
\end{subequations}
Additionally, for facilitating our subsequent analysis, it is beneficial to add (\ref{alg:Bregman_iter_regular_dops_fixednet_forbacksplit_equi_seq_x}) with (\ref{alg:Bregman_iter_regular_dops_fixednet_forbacksplit_equi_seq_y}), which yields:
\begin{subequations}\label{alg:Bregman_iter_regular_dops_fixednet_forbacksplit_equi_seq_simplified}
\begin{align}
\gamma  y_{k+1}-W(x_{k+1}-x_k)\in\gamma\partial f(x_{k+1}), \label{alg:Bregman_iter_regular_dops_fixednet_forbacksplit_equi_seq_simplified_x}\\
(I-W) x_{k+1}+\gamma(y_{k+1}-y_k)&= 0. \label{alg:Bregman_iter_regular_dops_fixednet_forbacksplit_equi_seq_simplified_y}
\end{align}
\end{subequations} 
It is clear from the above equation that the proposed algorithm after splitting is no longer the exact Bregman iterative method since $y_{k+1}\notin\partial f(x_{k+1})$. Thus, it also can be understood as an inexact version of Bregman iterative regularization.

We summarize the above proposed algorithm--Forward-Backward Bregman Splitting (D-FBBS)--in \emph{Algorithm}~\ref{alg:FBBS4fixednet}.
\begin{algorithm}[H]
\caption{D-FBBS for Fixed Networks}\label{alg:splitBregman_fixednet_detail}
\begin{algorithmic}[1]
\State {\bf Initialization}: $y_{i,0}=0, \forall i\in\mathcal{V}$ such that $\ones^Ty_0=0$, while the initial guess of $x_0$ can be arbitrarily assigned.
\State {\bf Primal Update}: For each agent $i\in\mathcal{V}$, compute:
\begin{equation*}
\begin{aligned}
x^{av}_{i,k}&=\sum_{j\in\mathcal{N}_i\cup \{i\}} w_{ij}x_{j,k},\\
x_{i,k+1}&=\underset{x_i\in\mathcal{R}^{d}}{\text{arg}\min}~\bracket{D^{y_{i,k}}_{f_i}(x_i,x_{i,k})+\frac{1}{2\gamma}\norm{x_i-x^{av}_{i,k}}^2}.
\end{aligned}
\end{equation*}
\State {\bf Dual Update}: For each agent $i\in\mathcal{V}$, $$y_{i,{k+1}}=y_{i,k}-\frac{1}{\gamma}\sum_{j\in\mathcal{N}_i}w_{ij}\bracket{x_{i,k+1}-x_{j,k+1}}.$$
\State Set k $\rightarrow$ k+1 and go to Step 2 until certain stopping criteria (e.g., maximum number of iterations) is satisfied.
\end{algorithmic}
\label{alg:FBBS4fixednet}
\end{algorithm}
\begin{rem}
\alert{It is obvious that the above algorithm can be carried out in a distributed manner since each agent only requires local information to solve its own optimization subproblem at each iteration. In particular, at time $k$, each agent $i\in\mathcal{V}$ keeps a local variable $x_{i,k}$ to be shared with its neighbors. Once each agent $i$ collects the information $x_{j,k}$ from its neighbors $j\in\mathcal{N}_i$, it solves the local optimization subproblem~(\ref{alg:Bregman_iter_regular_dops_fixednet_forbacksplit_x}) based on the obtained weighted average $x_{i,k}^{av}$, yielding the optimal solution $x_{i,k+1}$ which is then to be communicated to its neighbors. Upon receiving the estimated optimum $x_{j,k+1}$ from its neighbors, each agent $i$ updates the dual variable $y_{i,k}$ according to the step~(\ref{alg:Bregman_iter_regular_dops_fixednet_forbacksplit_y}) which can also be done locally.}
\end{rem}
\subsection{Some Basic Lemmas}
The following lemmas will be used in the subsequent analysis of the connections of the proposed algorithm to the existing ones as well as the theoretical proof of its convergence.
\begin{lem}\label{lem:bijective_transform}
Let $P$ be a $m\times m$ matrix such that $null(P)= span\{\ones\}$. Then, for each ${y}\in span^\perp\{\ones\}$, there exists a unique ${y}'\in span^\perp\{\ones\}$ such that ${y}=P{y}'$  and vice versa, i.e., the $P$-transformation between ${y}$ and ${y}'$ is bijective.
\end{lem}
\begin{IEEEproof}
See Appendix~\ref{sec:appendix}.
\end{IEEEproof}
\begin{lem}[Conservation Property]\label{lem:ConservProp}
Consider the sequence $\seq{y_k}$ generated by the algorithm~(\ref{alg:Bregman_iter_regular_dops_fixednet_forbacksplit}). Suppose $\ones^Ty_0=0$ and Assumption~\ref{assum:fixednet_weightmatrix} holds, then we have $\ones^Ty_k=0, \forall k\geq 0$.
\end{lem} 

The above lemma of conservation property is immediately followed by  multiplying $\ones^T$ from both sides of (\ref{alg:Bregman_iter_regular_dops_fixednet_forbacksplit_equi_seq_simplified_y}) and knowing from Assumption~\ref{assum:fixednet_weightmatrix} that $\ones^T(I-W)=0$.
\subsection{Variations and Connections to Existing Algorithms} 
\label{sec:variations_connections}
In this section, we first provide a variant of the proposed algorithm, termed \emph{Inexact} D-FBBS (ID-FBBS), to tackle cost functions having Lipschitz gradients and show their close connections with some existing well-known algorithms. In particular, under the proposed framework, we can perform another proper forward-backward splitting for functions having Lipschitz gradients, which essentially belongs to inexact Uzawa methods. We also show that, via proper change of variables, the proposed algorithm is equivalent to some existing well-known methods, such as preconditioned augmented Lagrangian methods and primal-dual approaches. 


\subsubsection*{\bf Inexact Uzawa Method} 
The proposed algorithm is still expensive in the sense that at each iteration we need to evaluate the inverse of $I+\gamma\partial f$. However, if we know that $f\in\Gamma_0(\mathcal{H})$ has Lipschitz gradients in the domain $\mathcal{D}$ of $f$,  i.e., 
\begin{equation}\label{eq:cost_lipshitzgradient}
\norm{\nabla f(x)-\nabla f(y)}\leq L_f\norm{x-y},~\forall x,y\in\mathcal{D},
\end{equation}
 then $f$ can be also evaluated explicitly in a cheaper way. 
That is, applying the forward-backward splitting technique gives
\begin{equation}\label{eq:splitBregman_dops_fixednet_smooth_x_forbacksplit_inclusion}
\begin{aligned}
x_{k+1}=Wx_k-\gamma(\nabla f(x_k)-y_k),
\end{aligned}
\end{equation} 
 which amounts to an inexact Uzawa Method being applied to the augmented Lagrangian~(\ref{eq:dops_lagrangian}), i.e., 
\begin{equation}\label{alg:inexact_Uzawa_to_Lagrangian}
\begin{aligned}
x_{k+1}&=x_k-\gamma\nabla_x \psi(x).\\
\end{aligned}
\end{equation} 
Thus, the above analysis leads to the following variant (\emph{Inexact} D-FBBS) of the proposed distributed algorithm:
\begin{subequations}\label{alg:dops_fixednet_ineaxct_Uzawa_smooth_forbacksplit}
\begin{align}
x_{k+1}&=Wx_k-\gamma(\nabla f(x_k)-y_k), \label{alg:dops_fixednet_ineaxct_Uzawa_smooth_forbacksplit_x}\\
y_{k+1}&=y_k-\frac{1}{\gamma} (I-W)x_{k+1}.  \label{alg:dops_fixednet_ineaxct_Uzawa_smooth_forbacksplit_y}
\end{align}
\end{subequations} 
Suppose $y_0$ is set to $0$. Then, summing~(\ref{alg:dops_fixednet_ineaxct_Uzawa_smooth_forbacksplit_y}) over $k$ and substituting it into~(\ref{alg:dops_fixednet_ineaxct_Uzawa_smooth_forbacksplit_x}) yields $$\underbrace{x_{k+1}=Wx_k-\gamma\nabla f(x_k)}_{\text{DSM}}-\underbrace{\sum^{k}_{i=0}(I-W)x_i}_{\text{Correction}},$$ which can be termed corrected DSM and is, indeed, equivalent to EXTRA in terms of $x$-update with $W=\tilde{W}=\frac{I+W'}{2}$, where $\tilde{W}$ and $W'$ are two properly designed weight matrices ~\cite[Eq.~(2.13)]{shi2015extra}. Thus, the convergence result\footnote{Note that we can obtain the linear convergent result as well if the cost function is known to be also strongly convex.} in terms of $x$-update for the OCP problem follows from the similar analysis therein, i.e., $\gamma\leq\frac{2\lambda_{min}(W)}{L_f}$. However, under our Bregman splitting framework, we are able to further improve the rate to $o(1/k)$ for ID-FBBS, thus also for EXTRA (cf.~Prop.~\ref{prop:convergence_IDFBBS}).

\begin{algorithm}[H]
\caption{ID-FBBS for Fixed Networks}\label{alg:splitBregman_fixednet_detail_inexact}
\begin{algorithmic}[1]
\State {\bf Initialization}: $y_{i,0}=0, \forall i\in\mathcal{V}$ such that $\ones^Ty_0=0$, while the initial guess of $x_0$ can be arbitrarily assigned.
\State {\bf Primal Update}: For each agent $i\in\mathcal{V}$, compute:
\begin{equation*}
\begin{aligned}
x_{i,k+1}=\sum_{j\in\mathcal{N}_i\cup \{i\}} w_{ij}x_{j,k}-\gamma(\nabla f_i(x_{i,k})-y_{i,k}).
\end{aligned}
\end{equation*}
\State {\bf Dual Update}: For each agent $i\in\mathcal{V}$, $$y_{i,{k+1}}=y_{i,k}-\frac{1}{\gamma}\sum_{j\in\mathcal{N}_i}w_{ij}\bracket{x_{i,k+1}-x_{j,k+1}}.$$
\State Set k $\rightarrow$ k+1 and go to Step 2 until certain stopping criteria (e.g., maximum number of iterations) is satisfied.
\end{algorithmic}
\label{alg:FBBS4fixednet_inexact}
\end{algorithm}


Now, we show that the proposed algorithm is, in fact, equivalent to some well-known existing algorithms by using preconditioned technique or proper change of variables.

\subsubsection*{\bf Preconditioned Augmented Lagrangian Method} To avoid the computation of the inverse of the weight matrix $W$, a clever way is to add an extra proximity term as follows
\begin{subequations}\label{eq:equl_precond_augLagrangian}
\begin{align}
x_{k+1}&=\underset{x\in\mathcal{R}^{m}}{\text{arg}\min}~\psi(x,y_k)+\frac{1}{2\gamma}\norm{x-x_k}^2_{W} \label{alg:dops_fixednet_precond_augLagrangian_forbacksplit_x}\\
y_{k+1}&=y_k-\frac{1}{\gamma}(I-W) x_{k+1}, \label{alg:dops_fixednet_precond_augLagrangian_forbacksplit_y}
\end{align}
\end{subequations}
where $\psi$ is the Lagrangian defined in~(\ref{eq:dops_lagrangian}). Introducing the prox-term allows the $x$-update step to be evaluated explicitly. Note that in the algorithm of Bregman Operator Splitting the above prox-term is replaced with the Bregman distance induced by a strongly convex function~\cite{zhang2010bregmanized,zhang2011unified}.  It is not difficult to verify that the above algorithm with $\ones^Ty_0=0$ is equivalent to the proposed distributed \emph{Algorithm}~\ref{alg:splitBregman_fixednet_detail}. Moreover, setting $y_0=0$, summing~(\ref{alg:dops_fixednet_precond_augLagrangian_forbacksplit_y}) over $k$ and substituting it into~(\ref{alg:dops_fixednet_precond_augLagrangian_forbacksplit_x}) yields $$x_{k+1}={\bf prox}_{\gamma f}\bracket{Wx_k-\sum^{k}_{i=0}(I-W)x_i},$$ which, similar as before, can be shown to be equivalent to P-EXTRA with $W=\tilde{W}=\frac{I+W'}{2}$ in terms of $x$-update~\cite{shi2015proximal}. However, we will show that our algorithm not only solve the OCP (primal) problem  but also the OEP (dual) problem. Thus, our convergence analysis is also different from theirs.

 Moreover, it is not difficult to verify that the above algorithm is also equivalent to the Jacobi variant of distributed AL methods proposed in~\cite{jakovetic2015linear} when the inner iterations of consensus are carried out only once, i.e., $\tau=1$. Note that their convergence result does not support this case.
\subsubsection*{\bf Primal-Dual Approach} 
In~\cite{chambolle2011first}, a general \emph{primal-dual} proximal point algorithm is proposed as follows
\begin{equation}\label{eq:primal_dual_approach}
\begin{aligned}
x_{k+1}&={\bf prox}_{\tau f}(x_k+\tau K\check{y}_k)\\
y_{k+1}&={\bf prox}_{\delta g}(y_k-\delta K^T x_{k+1})\\
\check{y}_{k+1}&=y_{k+1}+\theta (y_{k+1}-y_k),~\theta\in[0,1]
\end{aligned}
\end{equation}
to solve some generic saddle point problems encountered in image processing. In particular, when $g=0$ and $\theta=1$, it solves the following \emph{symmetric} saddle point problem:
\begin{equation}\label{eq:symmetric_saddlepoint_probs}
0\in\Bracket{
\begin{bmatrix}
\partial f &K\\
K^T &O
\end{bmatrix}\begin{bmatrix} x\\ y \end{bmatrix}}.
\end{equation}
Let $K=\sqrt{I-W}$. Since $null\{K\}=null\{I-W\}=span\{\ones\}$, we know from Lemma~\ref{lem:bijective_transform} that for each $y'\in span^\perp\{\ones\}$ there exists a unique $y\in span^\perp\{\ones\}$ such that $y'=Ky$. Thus, the \emph{primal-dual} algorithm~(\ref{eq:primal_dual_approach}) can be rewritten as follows \cite{condat2013primal,bianchi2014primal}:
\begin{equation}\label{eq:primal_dual_approach_singleobjective}
\begin{aligned}
x_{k+1}&={\bf prox}_{\tau f}(x_k+\tau (2y'_k-y'_{k-1})),\\
y'_{k+1}&=y'_k-\delta (I-W)x_{k+1}.\\
\end{aligned}
\end{equation}
It is easy to verify that when $\tau=\frac{1}{\delta}=\gamma, \ones^Ty_0=0$ (thus $\ones^Ty_k=0,\forall k$ by Lemma~\ref{lem:ConservProp}) the above algorithm is equivalent to the proposed distributed \emph{Algorithm}~\ref{alg:splitBregman_fixednet_detail} and it, in fact, solves the following \emph{asymmetric} saddle point problem:
\begin{equation}\label{eq:asymmetric_saddlepoint_probs}
0\in\Bracket{
\begin{bmatrix}
\partial f &I\\
I-W &O
\end{bmatrix}\begin{bmatrix} x\\ y' \end{bmatrix}},~\ones^Ty'=0,
\end{equation}
which is exactly the optimality conditions~(\ref{eq:dops_optimality_conditions}).
\begin{rem}
Note that the algorithm developed based on the symmetric formulation~(\ref{eq:symmetric_saddlepoint_probs}) is not suitable for stochastic networks as the saddle point varies with the topology, i.e., depending on the weight matrix $W$ and thus $K$. As such, the above connection via proper change of variables will be no longer valid for stochastic networks .
\end{rem}
\subsubsection*{\bf Decentralized ADMM}
\alert{The proposed algorithm is also related to decentralized ADMM (e.g., C-ADMM~\cite{mateos2010distributed} or D-ADMM~\cite{shi2014linear}) and their inexact versions (e.g., IC-ADMM~\cite{chang2015multi} or DLM~\cite{ling2015dlm}). The specific algorithm of D-ADMM/C-ADMM can be depicted in a compact form as follows:
\begin{equation}\label{alg:dadmm}
\begin{aligned}
x_{k+1}&=(2cD+\partial f)^{-1}\bracket{c(D+A)x_{k}-y_{k}}\\
y_{k+1}&=y_{k}+c(D-A)x_{k+1},
\end{aligned}
\end{equation}
where $A$ is the adjacency matrix and $D$ is the diagonal (degree) matrix corresponding to the Laplacian matrix $L$ of the graph. If the matrix $D$ is invertible, it can be rewritten as
\begin{equation}
\begin{aligned}
x_{k+1}&=(I+\frac{1}{2c}D^{-1}\partial f)^{-1}\bracket{(I-\frac{1}{2}L^{rw})x_{k}-\frac{1}{2c}y'_{k}}\\
y'_{k+1}&=y'_{k}+c(I+W^{rw})x_{k+1},
\end{aligned}
\end{equation}
where $L^{rw}=D^{-1}L$ is the random walk normalized Laplacian matrix, $W^{rw}=I-L^{rw}$ the corresponding weight matrix and $y'=D^{-1}y$ the scaled new dual variable. From the above analysis, it follows that the decentralized ADMM has a very similar form with our algorithm except for both the proximal parameter and the Laplacian matrix being scaled by $D^{-1}$.
Additionally, the inexact decentralized ADMM, such as DLM/IC-ADMM, can be generally depicted\footnote{Without loss of generality, for the IC-ADMM algorithm, we assume all nodes are using the same value of the parameter $\beta$ for update~\cite{chang2015multi}.} as follows:
\begin{equation}\label{alg:icadmm_dlm}
\begin{aligned}
x_{k+1}&=x_{k}-\tilde{D}^{-1}\bracket{\nabla f(x_{k})+c(D-A)x_{k}+y_{k}}\\
y_{k+1}&=y_{k}+c(D-A)x_{k+1},
\end{aligned}
\end{equation}
where $\tilde{D}=2cD+\rho I$ with $\rho$ being a parameter. Using the similar analysis as above, we can show that these inexact algorithms has a very similar form with ID-FBBS, the inexact version of D-FBBS, up to a scaling of $D^{-1}$ as shown before. 
\begin{rem}
Alternatively, if we replace the consensus constraint $(I-W)x=0$ with $(D-A)x=0$ in the formulation of the OCP problem, we can then easily recover all the above algorithms using the proposed Bregman splitting framework.. Additionally, the diagonal matrix $D$ has great impacts on the convergence performance and the parameter tuning for the above decentralized ADMM-based algorithms. In particular, when the degree of certain agent becomes zero (e.g., some agent being momentarily disconnected from the network), the algorithm will not be able to work properly. Instead, our algorithm, by effectively introducing the weight matrix, is more robust to the variation of the degree matrix (cf.~Section~\ref{sec:simulations}), implying that it is more suitable for varying networks.
\end{rem}}

\section{Distributed Optimization over Fixed Networks}
This section deals with distributed optimization problems over fixed networks in which the weight matrix is not changing with time. In particular, we study the convergence properties of the proposed algorithm for fixed networks.

Before proceeding to the main result, let us first establish the following lemmas:
\begin{lem}\label{lem:fixednet_fixedpointresidual_non_increasing}
Let $\tilde{x}=(I-\avector)x$, $u_{k+1}=\norm{\tilde{x}_{k+1}}^2_{I-W}+\norm{x_{k+1}-x_k}^2_W$ and suppose Assumptions~\ref{assum:fixednet_weightmatrix}~and~\ref{assum:costfuctions_proper_closed_convex} hold. Then, the sequence $\seq{u_k}$ generated by the D-FBBS algorithm (\ref{alg:Bregman_iter_regular_dops_fixednet_forbacksplit}) is monotonically non-increasing. We even have
\begin{equation}\label{eq:fixednet_iterseq_fixedpointresidual_U}
u_{k+1}\leq u_{k}-\norm{x_{k+1}-x_k}^2_{I-W}.
\end{equation}
\end{lem}

\begin{IEEEproof}
Let $\bar{x}\in span\{\ones\}$. Knowing that $\partial f$ is maximally monotone in the domain $\mathcal{R}^{m}$ by Assumption~\ref{assum:costfuctions_proper_closed_convex} and that $\gamma y_{k+1}-W(x_{k+1}-x_{k})\in\gamma\partial f(x_{k+1})$ and $\gamma y_k-W(x_k-x_{k-1})\in\gamma\partial f(x_k)$ from~(\ref{alg:Bregman_iter_regular_dops_fixednet_forbacksplit_equi_seq_simplified_x}),  we have
\begin{equation}\label{eq:fixednet_iterseq_monotone_f_kplus1_k}
\begin{aligned}
&\innprod{\gamma(y_{k+1}-y_k)-W(x_{k+1}-x_k)\right.\\
&~~~~~~~~~~~~~~~~~~~~~~~\left.+W(x_k-x_{k-1})}{x_{k+1}-x_k}\\
&=\innprod{-W(x_{k+1}-x_k)+W(x_k-x_{k-1})}{x_{k+1}-x_k}\\
&~~~~~-\innprod{(I-W)(x_{k+1}-\bar{x})}{{x_{k+1}-x_k}}\geq 0,\\
\end{aligned}
\end{equation}
where we have used (\ref{alg:Bregman_iter_regular_dops_fixednet_forbacksplit_equi_seq_simplified_y}) and the stochasticity of $W$ by Assumption~\ref{assum:fixednet_weightmatrix} to obtain the last term.

Recalling that $W$ is symmetric and positive definite by Assumption~\ref{assum:fixednet_weightmatrix}, the above inequality further leads to
\begin{equation}\label{eq:fixednet_iterseq_monotone_f_kplus1_k_equi}
\begin{aligned}
&\innprod{(I-W)(x_{k+1}-\bar{x})}{{x_{k+1}-x_k}}\\
&\leq\innprod{-W(x_{k+1}-x_k)+W(x_k-x_{k-1})}{x_{k+1}-x_k}\\
&=-\norm{x_{k+1}-x_k}^2_W+\innprod{W(x_k-x_{k-1})}{x_{k+1}-x_k}\\
&\leq-\norm{x_{k+1}-x_k}^2_W+\frac{\norm{x_{k+1}-x_k}^2_W+\norm{x_k-x_{k-1}}^2_W}{2}\\
&=\frac{-\norm{x_{k+1}-x_k}^2_W+\norm{x_k-x_{k-1}}^2_W}{2}.
\end{aligned}
\end{equation}
Combining the following identity
\begin{equation}\label{eq:fixednet_iterseq_identity}
\begin{aligned}
&2\innprod{(I-W)(x_{k+1}-\bar{x})}{{x_{k+1}-x_k}}\\
&=\norm{x_{k+1}-\bar{x}}^2_{I-W}-\norm{x_k-\bar{x}}^2_{I-W}+\norm{x_{k+1}-x_k}^2_{I-W}\\
\end{aligned}
\end{equation}
with~(\ref{eq:fixednet_iterseq_monotone_f_kplus1_k_equi}) yields
\begin{equation}
\begin{aligned}
&\norm{x_{k+1}-\bar{x}}^2_{I-W}-\norm{x_k-\bar{x}}^2_{I-W}+\norm{x_{k+1}-x_k}^2_{I-W}\\
&~~~~~~~\leq -\norm{x_{k+1}-x_k}^2_W+\norm{x_k-x_{k-1}}^2_W,
\end{aligned}
\end{equation}
which is equivalent to
\begin{equation}\label{eq:fixednet_iterseq_fixedpointresidual}
\begin{aligned}
&\norm{\tilde{x}_{k+1}}^2_{I-W}-\norm{\tilde{x}_k}^2_{I-W}+\norm{x_{k+1}-x_k}^2_{I-W}\\
&~~~~~~~~~~~~~\leq -\norm{x_{k+1}-x_k}^2_W+\norm{x_k-x_{k-1}}^2_W.
\end{aligned}
\end{equation}
Let $u_{k+1}=\norm{\tilde{x}_{k+1}}^2_{I-W}+\norm{x_{k+1}-x_k}^2_W$. Rearranging the terms of (\ref{eq:fixednet_iterseq_fixedpointresidual}) leads to (\ref{eq:fixednet_iterseq_fixedpointresidual_U}).
\end{IEEEproof}

Now, we are ready to present the main convergence result for the D-FBBS algorithm over fixed networks.
\begin{thm}\label{thm:convergence_mainresult_fixednet}
Suppose Assumptions~\ref{assum:fixednet_weightmatrix}, \ref{assum:costfuctions_proper_closed_convex} and \ref{assum:existence_saddlepoint}  hold. Then, the sequence $\seq{(x_k,y_k)}$ generated by the proposed D-FBBS algorithm~(cf.~Alg.~\ref{alg:splitBregman_fixednet_detail}) will converge to a saddle point $(x^\star,y^\star)$ of the \emph{primal-dual} problem. Moreover, the fixed point residual in terms of $\norm{\tilde{x}_{k+1}}^2_{I-W}+\norm{x_{k+1}-x_k}^2_W$ will decease at a non-ergodic rate of $o(\frac{1}{k})$. 
\end{thm}

\begin{IEEEproof}
Let $\mathcal{S}$ be the optimal solution set of the \emph{primal-dual} problem and $(x^\star,y^\star)\in\mathcal{S}$ be a saddle point. Since $\ones^Ty_k=0, \forall k\geq 0$ by Lemma~\ref{lem:ConservProp} and $null(I-W)=span\{\ones\}$, there exists a unique $y'_k\in span^\perp\{\ones\}$ such that $y_k=(I-W)y'_k$ by Lemma~\ref{lem:bijective_transform}. Knowing that $\partial f$ is maximally monotone by Assumption~\ref{assum:costfuctions_proper_closed_convex} and $\gamma y^\star\in\gamma\partial f(x^\star)$ from the optimality conditions~(\ref{eq:dops_optimality_conditions}), together with~(\ref{alg:Bregman_iter_regular_dops_fixednet_forbacksplit_equi_seq_simplified_x}) we have
\begin{equation}
\begin{aligned}
&\innprod{\gamma(y_{k+1}-y^\star)-W(x_{k+1}-x_k)}{x_{k+1}-x^\star}\\
&=\innprod{\gamma(I-W)(y'_{k+1}-y'^\star)}{x_{k+1}-x^\star}\\
&~~~~~~~~~~~~~~-\innprod{W(x_{k+1}-x_k)}{x_{k+1}-x^\star}\geq 0.\\
\end{aligned}
\end{equation}
Since $W$ is symmetric by Assumption~\ref{assum:fixednet_weightmatrix}, with~(\ref{alg:Bregman_iter_regular_dops_fixednet_forbacksplit_equi_seq_simplified_y}) we have
\begin{equation}\label{eq:fixednet_mainresult_iterseq_fixedpointresidual_raw}
\begin{aligned}
&\innprod{\gamma(y_{k+1}-y^\star)-W(x_{k+1}-x_k)}{x_{k+1}-x^\star}\\
&=\gamma\innprod{(I-W)(x_{k+1}-x^\star)}{y'_{k+1}-y'^\star}\\
&~~~~~~~~~~~~~~-\innprod{W(x_{k+1}-x_k)}{x_{k+1}-x^\star}\\
&=\gamma^2\innprod{-(y_{k+1}-y_k)}{y'_{k+1}-y'^\star}\\
&~~~~~~~~~~~~~~-\innprod{W(x_{k+1}-x_k)}{x_{k+1}-x^\star}\\
&=\gamma^2\innprod{-(I-W)(y'_{k+1}-y'_k)}{y'_{k+1}-y'^\star}\\
&~~~~~~~~~~~~~~-\innprod{W(x_{k+1}-x_k)}{x_{k+1}-x^\star}\geq 0.\\
\end{aligned}
\end{equation}
Then, using the similar identity as~(\ref{eq:fixednet_iterseq_identity}) yields
\begin{equation}\label{eq:fixednet_mainresult_iterseq_fixedpointresidual_main}
\begin{aligned}
&\gamma^2\norm{y'_{k+1}-y'^\star}^2_{I-W}-\gamma^2\norm{y'_k-y'^\star}^2_{I-W}+\norm{x_{k+1}-x^\star}^2_W\\
&-\norm{x_k-x^\star}^2_W\leq -\gamma^2\norm{y'_{k+1}-y'_k}^2_{I-W}-\norm{x_{k+1}-x_k}^2_W.
\end{aligned}
\end{equation}
Letting $\tilde{x}_k=(I-\avector)x_k$ and recalling that $y_{k+1}=(I-W)y'_{k+1}$ and $y_k=(I-W)y'_k$, we have  
\begin{equation}\label{eq:fixednet_mainresult_coordinatetransform_disagreement}
\begin{aligned}
\gamma^2&\norm{y'_{k+1}-y'_k}^2_{I-W}=\gamma^2\innprod{y'_{k+1}-y'_k}{y_{k+1}-y_k}\\
&~~~~~~~~~~=-\gamma\innprod{y'_{k+1}-y'_k}{(I-W)x_{k+1}}\\
&~~~~~~~~~~=-\innprod{\gamma(I-W)(y'_{k+1}-y'_k)}{x_{k+1}}\\
&~~~~~~~~~~=\norm{x_{k+1}}^2_{I-W}=\norm{\tilde{x}_{k+1}}^2_{I-W}.
\end{aligned}
\end{equation}
where we have used (\ref{alg:Bregman_iter_regular_dops_fixednet_forbacksplit_equi_seq_simplified_y}) to obtain the last two equalities. \alert{In addition, since $\ones^Ty'_k=0, \forall k\geq 0$ and $\rho\left(W-\avector\right)<1$ by Assumption~\ref{assum:fixednet_weightmatrix}, we have $y_{k+1}-y^\star=(I-W)(y'_{k+1}-y'^\star)=[I-(W-\avector)](y'_{k+1}-y'^\star)$ and in turn $$\norm{y'_{k+1}-y'^\star}^2_{I-W}=\norm{y_{k+1}-y^\star}^2_{[I-(W-\avector)]^{-1}}.$$
Thus,  the above relation (\ref{eq:fixednet_mainresult_iterseq_fixedpointresidual_main}) can be rewritten as
\begin{equation}\label{eq:fixednet_mainresult_iterseq_fixedpointresidual_main_equi}
\begin{aligned}
&\gamma^2\norm{y_{k+1}-y^\star}^2_{[I-(W-\avector)]^{-1}}+\norm{x_{k+1}-x^\star}^2_W\\
&~~~~~-\gamma^2\norm{y_k-y^\star}^2_{[I-(W-\avector)]^{-1}}-\norm{x_k-x^\star}^2_W\\
&~~~~~~~~~~~~~~~\leq -\norm{\tilde{x}_{k+1}}^2_{I-W}-\norm{x_{k+1}-x_k}^2_W.
\end{aligned}
\end{equation}}
Let $V_k=\gamma^2\norm{y_k-y^\star}^2_{[I-(W-\avector)]^{-1}}+\norm{x_k-x^\star}^2_W$ be the Lyapunov function, which is positive by Assumption~\ref{assum:fixednet_weightmatrix}. 
Then, summing (\ref{eq:fixednet_mainresult_iterseq_fixedpointresidual_main_equi}) over $k$ from $0$ to $n-1$ yields
\begin{equation}
\begin{aligned}
&\sum_{k=0}^{n-1}\left(\norm{\tilde{x}_{k+1}}^2_{I-W}+\norm{x_{k+1}-x_k}^2_W\right)\leq V_0-V_{n}< \infty.
\end{aligned}
\end{equation}
By Lemma 2, $u_{k+1}=\norm{\tilde{x}_{k+1}}^2_{I-W}+\norm{x_{k+1}-x_k}^2_W$ is monotonically non-increasing. Thus, we have
\begin{equation}
\begin{aligned}
&nu_n\leq\sum_{k=0}^{n-1}u_{k+1}\leq V_0,
\end{aligned}
\end{equation}
yielding $u_n\leq \frac{V_0}{n}=o(\frac{1}{n})$.
\alert{Since $u_n\geq 0$ by Assumption~\ref{assum:fixednet_weightmatrix}, we claim that $\lim_{n\rightarrow\infty}u_n=0$, which further implies that $\lim_{k\rightarrow\infty}\bracket{x_{k+1}-x_k}=0$ and $\lim_{k\rightarrow\infty}\tilde{x}_{k+1}=0$. From (\ref{eq:fixednet_mainresult_iterseq_fixedpointresidual_main_equi}) we know that $V_k$ is bounded and so is $x_k$ and $y_k$. Thus, the sequence $\seq{(x_k,y_k)}$ has at least a cluster point, say $(x_\infty,y_\infty)$ and the corresponding subsequence $\{(x_{k_i},y_{k_i})\}_{i\geq 0}$. Taking the limit along $k_i\rightarrow\infty$ in (\ref{alg:Bregman_iter_regular_dops_fixednet_forbacksplit_equi_seq_simplified}) we immediately obtain that the cluster point $(x_\infty,y_\infty)$ is a saddle point of the \emph{primal-dual} problem (cf.~Prop.~\ref{prop:optim_iff_saddlepoint}). Since $\seq{V_k}$ converges and is contractive with respect to the optimal set $\mathcal{S}$, (\ref{eq:fixednet_mainresult_iterseq_fixedpointresidual_main_equi}) implies that $(x_\infty,y_\infty)$ is the unique cluster point (cf. \cite[Th.~5.5]{bauschke2011convex} and \cite[Th.~3.7]{he2012convergence}).  It then follows that $\seq{(x_k,y_k)}$ will converge to a unique point of the optimal set $\mathcal{S}$. Thus, we conclude that the sequence $\seq{(x_k,y_k)}$ will converge to a saddle point $(x^\star,y^\star)$  with a non-ergodic rate of $o(\frac{1}{k})$.}
\end{IEEEproof}

The following proposition is to show the convergence of the ID-FBBS algorithm over fixed networks.
\begin{prop}\label{prop:convergence_IDFBBS}
\alert{Suppose Assumptions~\ref{assum:fixednet_weightmatrix}, \ref{assum:costfuctions_proper_closed_convex}, \ref{assum:existence_saddlepoint}  hold, and the cost function $f$ has $L_f$-Lipschitz gradient (cf.~Eq.~(\ref{eq:cost_lipshitzgradient})). Let $\gamma<\frac{\lambda_{\min}(W)}{L_f}$. Then, the sequence $\seq{(x_k,y_k)}$ generated by the proposed ID-FBBS algorithm (cf. Alg.~\ref{alg:splitBregman_fixednet_detail_inexact}) will converge to a saddle point $(x^\star,y^\star)$ of the \emph{primal-dual} problem. Moreover, the fixed point residual in terms of $\norm{\tilde{x}_{k}}^2_{I-W}+\norm{x_{k+1}-x_k}^2_W$will decrease at a non-ergodic rate of $o(\frac{1}{k})$. }
\end{prop}
\begin{IEEEproof}
See Appendix~\ref{sec:appendix}.
\end{IEEEproof}

\begin{rem}
\alert{ID-FBBS is dedicated for cost functions that have Lipschitz gradients. It can be regarded as a gradient-based method and is thus more computationally efficient compared to D-FBBS which is designed for general convex cost functions. Both algorithms have a non-ergodic convergence rate of $o(1/k)$, the best known rate in the existing literature.}
\end{rem}


\section{Distributed Optimization over Stochastic Networks}
This section deals with the application of the proposed algorithm to stochastic networks. In particular, we show that under stronger assumption on the cost functions, i.e., being strongly convex, the same proposed distributed algorithm can be employed to solve
 the distributed optimization problem over stochastic networks.  Note that most of the analysis for a fixed network is not readily transferable to stochastic networks as the Lyapunov function employed therein is dependent on the network (i.e., the weight matrix $W$) thus varying with time. We have to find a new (common) Lyapunov function that is immune to varying networks. Bregman distance will thus play a key role in the sequent convergence analysis. 

We let $W_k$ be the weight matrix employed for communication by agents at time $k$. Each entry of $W_k$ is positive only if there is an active communication link between the associated agents at time $k$. Similar with fixed networks, we make the following assumptions.
\begin{assum}\label{assum:asynet_weightmatrix}
Let $\seq{W_k}$ be the i.i.d. stochastic weight matrix sequence and $\overbar{W}=\expect{W_k}$ be the mean. In addition, the following conditions hold: $W_k^T=W_k,~W_k>0,~W_k\ones=\ones,~\forall k\geq0,~\text{and}~\rho\left(\overbar{W}-\avector\right)<1$. 
\end{assum}
\begin{assum}\label{assum:asynet_costfunctions_stronglymonotone}
The cost functions are strongly convex such that $D^{\partial f(y)}_{f}(x,y)\geq \frac{\alpha}{2}\norm{x-y}^2$.
\end{assum}

We use the same algorithm as~(\ref{alg:Bregman_iter_regular_dops_fixednet_forbacksplit}) to solve the following stochastic optimal consensus problem (SOCP):
\begin{equation*}\label{prob:dops_asynet}
\begin{aligned}
(SOCP)~~~\min_{x\in\mathcal{R}^{md}} f(x)=\sum_{i=1}^mf_i(x_i) ~~~~~{\bf s.t.}~(I-\overbar{W})x=0.
\end{aligned}
\end{equation*}
The specific algorithm is as follows:
\begin{equation}\label{alg:Bregman_iter_regular_dops_asynet_forbacksplit}
\begin{aligned}
x_{k+1}&=\underset{x\in\mathcal{R}^{m}}{\text{arg}\min}~\bracket{D^{y_k}_f(x,x_k)+\frac{1}{2\gamma}\norm{x-W_kx_k}^2}\\
y_{k+1}&=y_k-\frac{1}{\gamma}(I-W_k) x_{k+1},\\
\end{aligned}
\end{equation}
which is equivalent to
\begin{subequations}\label{alg:Bregman_iter_regular_dops_asynet_forbacksplit_equi}
\begin{align}
\gamma y_k-(x_{k+1}-W_kx_k)&\in\gamma\partial f(x_{k+1}), \label{alg:Bregman_iter_regular_dops_asynet_forbacksplit_equi_seq_x}\\
(I-W_k) x_{k+1}+\gamma(y_{k+1}-y_k)&= 0. \label{alg:Bregman_iter_regular_dops_asynet_forbacksplit_equi_seq_y}
\end{align}
\end{subequations} 
As before, adding (\ref{alg:Bregman_iter_regular_dops_asynet_forbacksplit_equi_seq_x}) with (\ref{alg:Bregman_iter_regular_dops_asynet_forbacksplit_equi_seq_y}) yields:
\begin{subequations}\label{alg:Bregman_iter_regular_dops_asynet_forbacksplit_equi_seq_simplified}
\begin{align}
\gamma  y_{k+1}-W_k(x_{k+1}-x_k)\in\gamma\partial f(x_{k+1}) \label{alg:Bregman_iter_regular_dops_asynet_forbacksplit_equi_seq_simplified_x},\\
(I-W_k) x_{k+1}+\gamma(y_{k+1}-y_k)&= 0. \label{alg:Bregman_iter_regular_dops_asynnet_forbacksplit_equi_seq_simplified_y}
\end{align}
\end{subequations} 
We summarize the above proposed D-FBBS algorithm for stochastic networks in \emph{Algorithm}~\ref{alg:FBBS4asynet} as follows:
\begin{algorithm}[H]
\caption{D-FBBS for Stochastic Networks}\label{alg:splitBregman_asynet_detail}
\begin{algorithmic}[1]
\State {\bf Initialization}: $y_{i,0}=0, \forall i\in\mathcal{V}$ such that $\ones^Ty_0=0$, while the initial guess of $x_0$ can be arbitrarily assigned.
\State {\bf Primal Update}: For each agent $i\in\mathcal{V}$, compute:
\begin{equation*}
\begin{aligned}
x^{av}_{i,k}&=\sum_{j\in\mathcal{N}_i\cup \{i\}} w_{ij,k}x_{j,k},\\
x_{i,k+1}&=\underset{x_i\in\mathcal{R}^{d}}{\text{arg}\min}~\bracket{D^{y_{i,k}}_{f_i}(x_i,x_{i,k})+\frac{1}{2\gamma}\norm{x_i-x^{av}_{i,k}}^2}.
\end{aligned}
\end{equation*}
\State {\bf Dual Update}: For each agent $i\in\mathcal{V}$, $$y_{i,{k+1}}=y_{i,k}-\frac{1}{\gamma}\sum_{j\in\mathcal{N}_i}w_{ij,k}\bracket{x_{i,k+1}-x_{j,k+1}}.$$
\State Set k $\rightarrow$ k+1 and go to Step 2 until certain stopping criteria (e.g., maximum number of iterations) is satisfied.
\end{algorithmic}
\label{alg:FBBS4asynet}
\end{algorithm}
Now, we present the main result for stochastic networks.
\begin{thm}\label{thm:convergence_mainresult_rndnet}
Let $\lambda_{max}: =\max\{\seq{\lambda(W_k)}\}\in(0,1]$ and $\lambda_{min}: =\min\{\seq{\lambda(W_k)}\}\in(0,1]$ be the maximum and minimum eigenvalue of $W_k$ for all $k\geq 0$ respectively. Suppose Assumptions~\ref{assum:existence_saddlepoint},~\ref{assum:asynet_weightmatrix} and~\ref{assum:asynet_costfunctions_stronglymonotone} hold, and 
$$\gamma> \frac{2(1-\mu)(1-\lambda_{min})}{\alpha\mu}+\frac{2\lambda_{max}}{\alpha},$$
where $\mu\in(0,1)$. Then, the sequence $\seq{x_k}$ generated by the stochastic D-FBBS algorithm (cf.~Alg.~\ref{alg:splitBregman_asynet_detail}) will converge almost surely to an optimal solution of the SOCP problem. Moreover, let $\hat{x}_{k}=\frac{1}{k}\sum^{k-1}_{t=0}x_t$ and $\hat{\bar{x}}_k=\frac{1}{k}\sum^{k-1}_{t=0}\bar{x}_t$ be the running average of $x_t$ and $\bar{x}_t$ respectively. Then, the expected fixed point residual in terms of $\expect{\norm{\hat{x}_{k+1}-\hat{x}_k}^2_{(\frac{\gamma \alpha}{2}-\lambda_{max}-\eta)I}+\frac{1-\mu}{1-\nu}\norm{\hat{x}_k-\hat{\bar{x}}_k}^2_{I-\overbar{W}}}$ with $\nu\in(0,1)$, will decrease at an ergodic rate of $O(\frac{1}{k})$.
\end{thm}
\begin{IEEEproof}
Let us first consider an interesting recursive relation as follows:
\begin{equation} \label{eq:Bregman_threeidentity}
\begin{aligned}
&D^{q_{k+1}}_{\gamma f}(x,x_{k+1})-D^{q_{k}}_{\gamma f}(x,x_k)+D^{q_{k}}_{\gamma f}(x_{k+1},x_k)\\
&=\gamma f(x)-\gamma f(x_{k+1})-\innprod{q_{k+1}}{x-x_{k+1}}-\gamma f(x)+\gamma f(x_{k})\\
&~~~+\innprod{q_{k}}{x-x_{k}}+\gamma f(x_{k+1})-\gamma f(x_k)\\
&~~~-\innprod{q_{k}}{x_{k+1}-x_k}=\innprod{q_{k+1}-q_{k}}{x_{k+1}-x}.
\end{aligned}
\end{equation}
Let $q_{k}=\gamma y_{k}-W_{k-1}(x_{k}-x_{k-1})$. Using~(\ref{alg:Bregman_iter_regular_dops_asynnet_forbacksplit_equi_seq_simplified_y}) we have
\begin{equation} \label{eq:asynnet_mainresult_subgradient_difference}
\begin{aligned}
&q_{k+1}-q_{k}\\
&=\gamma(y_{k+1}-y_{k})-W_k(x_{k+1}-x_k)+W_{k-1}(x_k-x_{k-1})\\
&=-(I-W_k)x_{k+1}-W_k(x_{k+1}-x_k)+W_{k-1}(x_k-x_{k-1}).\\
\end{aligned}
\end{equation}
Let $\bar{x}\in span\{\ones\}$. Combining~(\ref{eq:Bregman_threeidentity}) and~(\ref{eq:asynnet_mainresult_subgradient_difference}) yields
\begin{equation}\label{eq:asynnet_mainresult_iterseq_main}
\begin{aligned}
&D^{q_{k+1}}_{\gamma f}(\bar{x},x_{k+1})-D^{q_{k}}_{\gamma f}(\bar{x},x_k)+D^{q_{k}}_{\gamma f}(x_{k+1},x_k)\\
&=\innprod{-W_k(x_{k+1}-x_k)+W_{k-1}(x_k-x_{k-1})}{x_{k+1}-\bar{x}}\\
&~~~~~-\norm{x_{k+1}-\bar{x}}^2_{I-W_k}\\
&=-\norm{x_{k+1}-\bar{x}}^2_{I-W_k}-\innprod{W_k(x_{k+1}-x_k)}{x_{k+1}-\bar{x}}\\
&~~~~~+\innprod{W_{k-1}(x_k-x_{k-1})}{x_k-\bar{x}}\\
&~~~~~+\innprod{W_{k-1}(x_k-x_{k-1})}{x_{k+1}-x_k}\\
&\leq-\norm{x_{k+1}-\bar{x}}^2_{I-W_k}-\innprod{W_k(x_{k+1}-x_k)}{x_{k+1}-\bar{x}}\\
&~~~~~+\innprod{W_{k-1}(x_k-x_{k-1})}{x_k-\bar{x}}-\norm{x_{k+1}-x_k}^2_\frac{{W_{k}}}{2}\\
&~~~~~+\norm{x_k-x_{k-1}}^2_\frac{{W_{k-1}}}{2}+\norm{x_{k+1}-x_k}^2_\frac{{W_{k}+W_{k-1}}}{2},\\
\end{aligned}
\end{equation}
where we have used the following inequality:
\begin{equation}
\begin{aligned}
&\innprod{W_{k-1}(x_{k}-x_{k-1})}{x_{k+1}-x_k}\\
&~~\leq \norm{x_{k+1}-x_k}^2_{\frac{W_{k-1}}{2}}+\norm{x_{k}-x_{k-1}}^2_{\frac{W_{k-1}}{2}}\\
&~~=\norm{x_{k+1}-x_k}^2_{\frac{{W_{k}+W_{k-1}}}{2}-\frac{{W_{k}}}{2}}+\norm{x_k-x_{k-1}}^2_\frac{{W_{k-1}}}{2}\\
\end{aligned}
\end{equation}
to obtain the last relation. Now, let 
\begin{equation}
\begin{aligned}
V_{k+1}&=D^{q_{k+1}}_{\gamma f}(\bar{x},x_{k+1})+\innprod{W_k(x_{k+1}-x_k)}{x_{k+1}-\bar{x}}\\
&~~~~+\norm{x_{k+1}-x_k}^2_\frac{{W_{k}}}{2},
\end{aligned}
\end{equation}
which is positive if $\gamma> \frac{\lambda_{max}}{\alpha}$ since
\begin{equation}\label{eq:asynnet_mainresult_positive_lyapunov}
\begin{aligned}
V_{k+1}&=D^{q_{k+1}}_{\gamma f}(\bar{x},x_{k+1})-\norm{x_{k+1}-\bar{x}}^2_\frac{{W_{k}}}{2}+\norm{x'_{k+1}-\bar{x}}^2_\frac{{W_{k}}}{2}\\
&\geq\norm{x_{k+1}-\bar{x}}^2_ {\frac{{\gamma\alpha I-W_{k}}}{2}}\geq\norm{x_{k+1}-\bar{x}}^2_{\frac{\gamma\alpha-\lambda_{max}}{2}I}>0,
\end{aligned}
\end{equation}
where $x'_{k+1}=2x_{k+1}-x_k$ and we have used Assumption~\ref{assum:asynet_costfunctions_stronglymonotone} to obtain $D^{q_{k}}_{\gamma f}(x_{k+1},\bar{x})\geq \frac{\gamma \alpha}{2}\norm{x_{k+1}-\bar{x}}^2$.

Then, (\ref{eq:asynnet_mainresult_iterseq_main}) can be rewritten as
\begin{equation}\label{eq:asynnet_mainresult_iterseq_main_simplified}
\begin{aligned}
&D^{q_{k}}_{\gamma f}(x_{k+1},x_k)-\norm{x_{k+1}-x_k}^2_\frac{{W_{k}+W_{k-1}}}{2}\\
&~~~~~~~~~~~~~~~~~~~~~+\norm{x_{k+1}-\bar{x}}^2_{I-W_k}\leq V_{k}-V_{k+1}.
\end{aligned}
\end{equation}
Applying the basic inequality in $G$-space $\norm{a+b}^2_G\geq (1-\frac{1}{\mu})\norm{a}^2_G+(1-\mu)\norm{b}^2_G,~\forall \mu\geq 0$, we have
\begin{equation}
\begin{aligned}
&\norm{x_{k+1}-\bar{x}}^2_{I-W_k}\geq\\
&~~~~~~~(1-\frac{1}{\mu})\norm{x_{k+1}-x_k}^2_{I-W_k}+(1-\mu)\norm{x_k-\bar{x}}^2_{I-W_k},
\end{aligned}
\end{equation}
where we require $\mu\in(0,1)$. Thus, (\ref{eq:asynnet_mainresult_iterseq_main_simplified}) becomes
\begin{equation}\label{eq:asynnet_mainresult_iterseq_main_simplified_relaxed}
\begin{aligned}
&D^{q_{k}}_{\gamma f}(x_{k+1},x_k)-\norm{x_{k+1}-x_k}^2_{(\lambda_{max}+\eta) I}\\
&~~~~~~~~~~~~~+(1-\mu)\norm{x_k-\bar{x}}^2_{I-W_k}\leq V_k-V_{k+1},
\end{aligned}
\end{equation}
where $\eta=\frac{(1-\mu)(1-\lambda_{min})}{\mu}$.

Knowing that $W_k$ is independent of the past states, taking total expectation of both sides yields
\begin{equation}\label{eq:asynnet_mainresult_iterseq_expect}
\begin{aligned}
&\expect{\norm{x_{k+1}-x_k}^2_{(\frac{\gamma \alpha}{2}-\lambda_{max}-\eta)I}}\\
&~~~~~+(1-\mu)\expect{\norm{x_k-\bar{x}_k}^2_{I-\overbar{W}}}\leq \expect{V_k}-\expect{V_{k+1}},
\end{aligned}
\end{equation}
where we have used the fact that $\norm{x_k-\bar{x}}^2_{I-\overbar{W}}=\norm{\tilde{x}_k}^2_{I-\overbar{W}}$ (cf. Assumption~\ref{assum:asynet_weightmatrix}) to replace the second term.

Summing~(\ref{eq:asynnet_mainresult_iterseq_expect}) over $k$ from $0$ through $n-1$ leads to
\begin{equation}\label{eq:asynnet_mainresult_iterseq_expect_summation}
\begin{aligned}
&\sum^{n-1}_{k=0}\expect{\norm{x_{k+1}-x_k}^2_{(\frac{\gamma \alpha}{2}-\lambda_{max}-\eta)I}}\\
&~~~~~+\sum^{n-1}_{k=0}(1-\mu)\expect{\norm{x_k-\bar{x}_k}^2_{I-\overbar{W}}}\leq \expect{V_0}<\infty,
\end{aligned}
\end{equation}
Suppose $\gamma>\frac{2(\lambda_{max}+\eta)}{\alpha}$ and let $n\rightarrow\infty$. Then, by Markov's inequality~\cite{Book:Allan2005_probability}, using $(a+b)^2\leq 2a^2+2b^2$ we have for $\forall \varepsilon>0$
\begin{equation}
\begin{aligned}
&\sum^\infty_{k=0}P(\norm{x_{k+1}-x_k}_{(\frac{\gamma \alpha}{2}-\lambda_{max}-\eta)I}+\sqrt{(1-\mu)}\norm{\tilde{x}_k}_{I-\overbar{W}}>\varepsilon)\\
&\leq\frac{2\sum^{\infty}_{k=0}\expect{\norm{x_{k+1}-x_k}^2_{(\frac{\gamma \alpha}{2}-\lambda_{max}-\eta)I}}}{\varepsilon^2}\\
&~~~+\frac{2\sum^{\infty}_{k=0}(1-\mu)\expect{\norm{x_k-\bar{x}_k}^2_{I-\overbar{W}}}}{\varepsilon^2}\leq\frac{2\expect{V_0}}{\varepsilon^2}<\infty.
\end{aligned}
\end{equation}
Thus, by Borel-Cantelli Lemma and Prop. 1.2 \cite[p.206]{Book:Allan2005_probability}, we have $\lim_{k\rightarrow\infty} x_{k+1}=x_k$ and $\lim_{k\rightarrow\infty} x_k=\bar{x}_k$ with probability one. \alert{Further, since $\gamma> \frac{2(1-\mu)(1-\lambda_{min})}{\alpha\mu}+\frac{2\lambda_{max}}{\alpha}$, from~(\ref{eq:asynnet_mainresult_iterseq_main_simplified_relaxed}) we know that $V_k$ is bounded and so is $x_k$ by (\ref{eq:asynnet_mainresult_positive_lyapunov}). Recalling that $q_k=\gamma y_{k}-W_{k-1}(x_{k}-x_{k-1})$, from~(\ref{alg:Bregman_iter_regular_dops_asynet_forbacksplit_equi_seq_simplified_x}) we have $q_k\in\gamma\partial f(x_{k})~\forall k$. Thus, by \cite[Prop. 4.2.3]{Book:ConvexOptim} we know that $q_k$ is bounded and so is $y_k$. Let $(x_\infty,y_\infty)$ be one of the cluster points of the sequence $\seq{(x_k,y_k)}$. Passing to the limit along (\ref{alg:Bregman_iter_regular_dops_asynet_forbacksplit_equi_seq_simplified}) we obtain that $(x_\infty,y_\infty)$ is a saddle point of the \emph{primal-dual} problem. Then, by  Proposition~\ref{prop:optim_iff_saddlepoint}, the limit point $x_\infty\in span\{\ones\}$ is also an optimal solution of the OCP problem. Since $\seq{V_k}$ converges and is contractive with respect to the consensus value $\bar{x}\in span\{\ones\}$, (\ref{eq:asynnet_mainresult_iterseq_main_simplified_relaxed}) implies that $x_\infty$ is the unique cluster point (cf. \cite[Th.~5.5]{bauschke2011convex} and \cite[Th.~3.7]{he2012convergence}). Thus, we conclude that $\seq{x_k}$ will converge almost surely to an optimal solution of the SCOP problem.}

Moreover, let $\hat{x}_{n}=\frac{1}{n}\sum^{n-1}_{k=0}x_k$. Multiplying both sides of~(\ref{eq:asynnet_mainresult_iterseq_expect_summation}) by $\frac{1}{n}$ and using the Jensen's inequality yields
\begin{equation}
\begin{aligned}
&\expect{\norm{\hat{x}_{n+1}-\hat{x}_n+\frac{1}{n}(\hat{x}_{n+1}-x_0)}^2_{(\frac{\gamma \alpha}{2}-\lambda_{max}-\eta) I}}\\
&~~~~~~~~~~~~~~+(1-\mu)\expect{\norm{\hat{x}_n-\hat{\bar{x}}_n}^2_{I-\overbar{W}}}\leq \frac{\expect{V_0}}{n},
\end{aligned}
\end{equation}
where $\hat{\bar{x}}_n=\frac{1}{n}\sum^{n-1}_{k=0}\bar{x}_k$. Again, using the basic inequality $\norm{a+b}^2\geq (1-\frac{1}{\nu})\norm{a}^2+(1-\nu)\norm{b}^2,~\forall \nu\geq 0$, we have
\begin{equation}
\begin{aligned}
&\expect{(1-\nu)\norm{\hat{x}_{n+1}-\hat{x}_n}^2_{(\frac{\gamma \alpha}{2}-\lambda_{max}-\eta) I}}\\
&~~~~~~~~~~~~~~~~~~~~~~~~+(1-\mu)\expect{\norm{\hat{x}_n-\hat{\bar{x}}_n}^2_{I-\overbar{W}}} \\
&\leq \frac{\expect{V_0}}{n}+\frac{1}{n^2}(\frac{1}{\nu}-1)\norm{\hat{x}_{n+1}-x_0}^2_{(\frac{\gamma \alpha}{2}-\lambda_{max}-\eta) I},
\end{aligned}
\end{equation}
where we require $\nu\in(0,1)$. Since $x_n$ is a convergent sequence and so is $\hat{x}_n$, thus $\hat{x}_n$ is uniformly bounded, implying that the last term of the above relation is of $O(\frac{1}{n^2})$. It follows that the expected fixed point residual in terms of the running average $\hat{x}_n$ will decrease at an ergodic rate of $O(\frac{1}{n})$.
\end{IEEEproof}




\section{Applications to Sensor Fusion}\label{sec:simulations}
We provide some simulations for fixed network and stochastic network respectively to show the effectiveness of the proposed algorithm.  In particular, we consider a canonical distributed estimation problem. Each sensor is assumed to measure certain unknown parameter $\theta\in\mathcal{R}^d$ with some Gaussian noise $\omega_i$, i.e., $z_i=M_i\theta+\omega_i$, where $M_i\in R^{s\times d}$ is the measurement matrix of sensor $i$ and $z_i\in R^s$ is the measurement data collected by sensor $i$. Thus, the maximum likelihood estimation with regularization can be casted as the following minimization problem:
\begin{equation}
\begin{aligned}
{\theta}^\star=\text{arg}\underset{\theta\in\mathcal{R}^d}{\min}~\bracket{\sum^m_{i=1}\norm{z_i-M_i\theta}^2+\lambda_{reg}\norm{\theta}^2},\\
\end{aligned}
\end{equation} 
where $\lambda_{reg}$ is the regularization parameter.  

\begin{figure}[!hbpt] 
      \centering
      \includegraphics[scale=0.54]{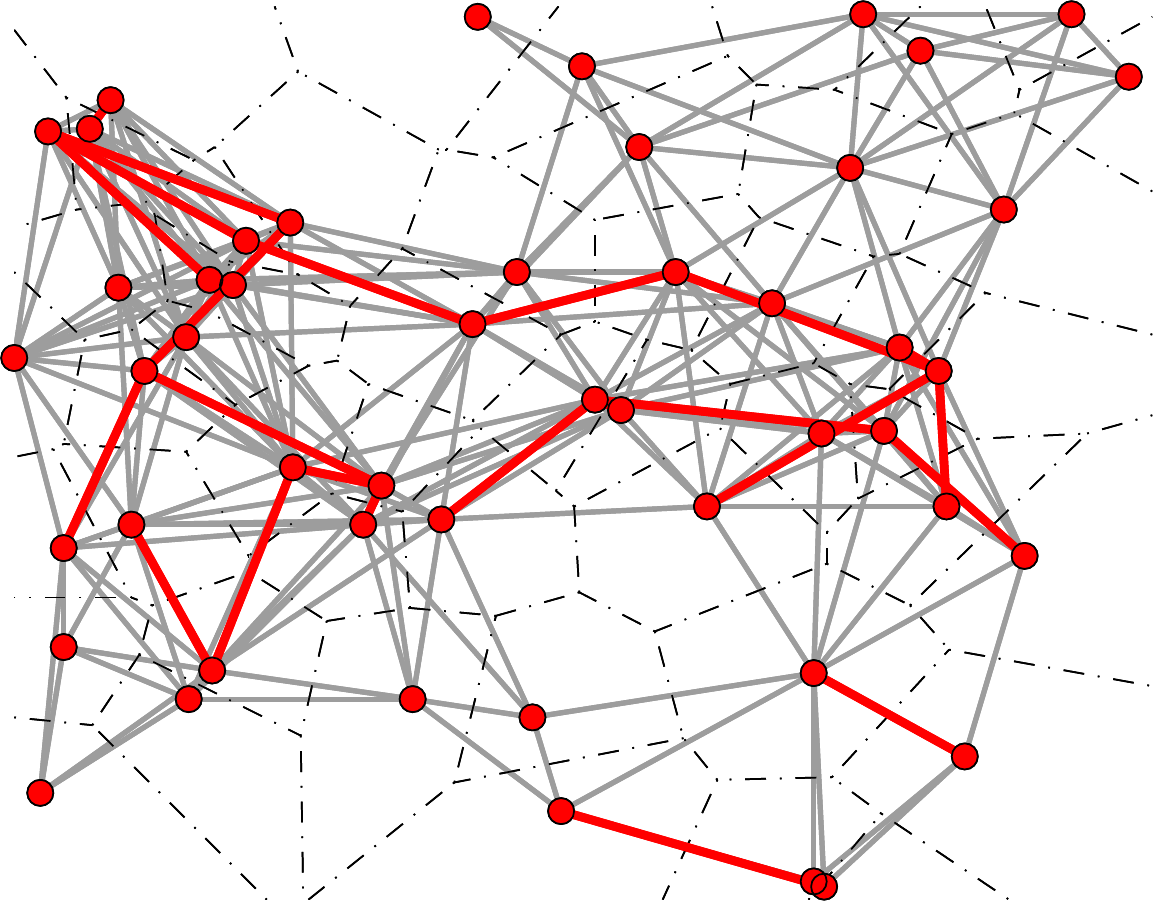}   
      \caption{A snapshot of a random sensor network of 50 nodes at certain time. The red lines denote the communication links being activated while the gray lines stand for no communication being carried out at this moment.}
      \label{fig:rndnet_snapshot}
\end{figure}

\alert{{\bf Network Setting:} We consider a random network whose nodes are randomly generated in a unit square and are connected with each other if they fall in the sensing range of nodes, i.e., $r\cdot\sqrt{\frac{\log{m}}{m}}$ (cf.~Fig.~\ref{fig:rndnet_snapshot}). The bigger the value of $r$, the more connected the graph will be and the larger of the maximum degree of the graph is. In addition, each communication link is assumed to be subject to random failures following certain Bernoulli Process. That is, at each iteration, each communication link will be activated with probability of $p$ and deactivated with $1-p$. Thus, when $p=1$, the random network will reduce to a fixed network.}

\begin{figure*}[!htpb]
\centering
 \subfloat[$r=1$]{ \includegraphics[width=0.32\textwidth]{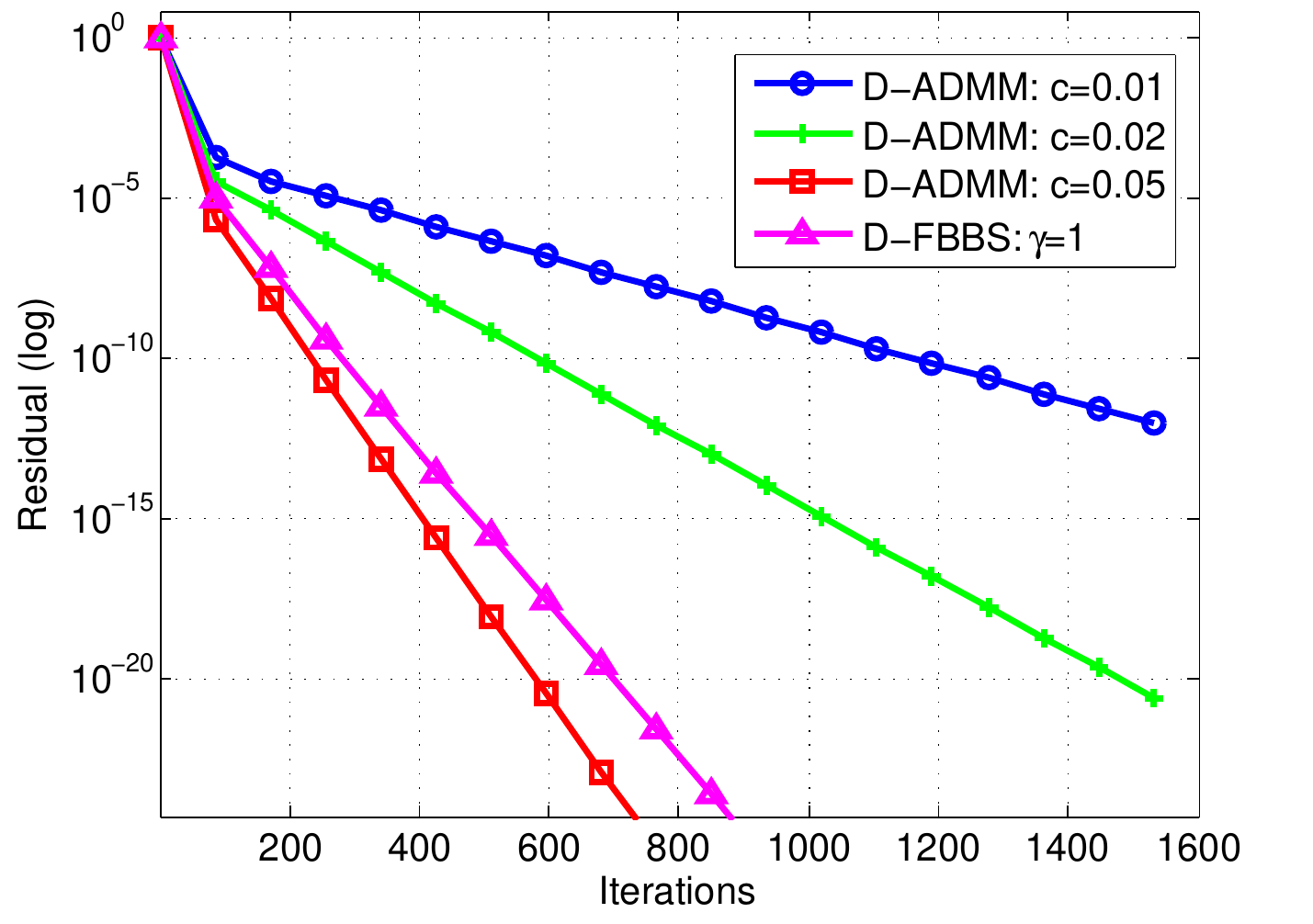}\label{fig:nocomm}}
 \subfloat[$r=1.5$]{ \includegraphics[width=0.32\textwidth]{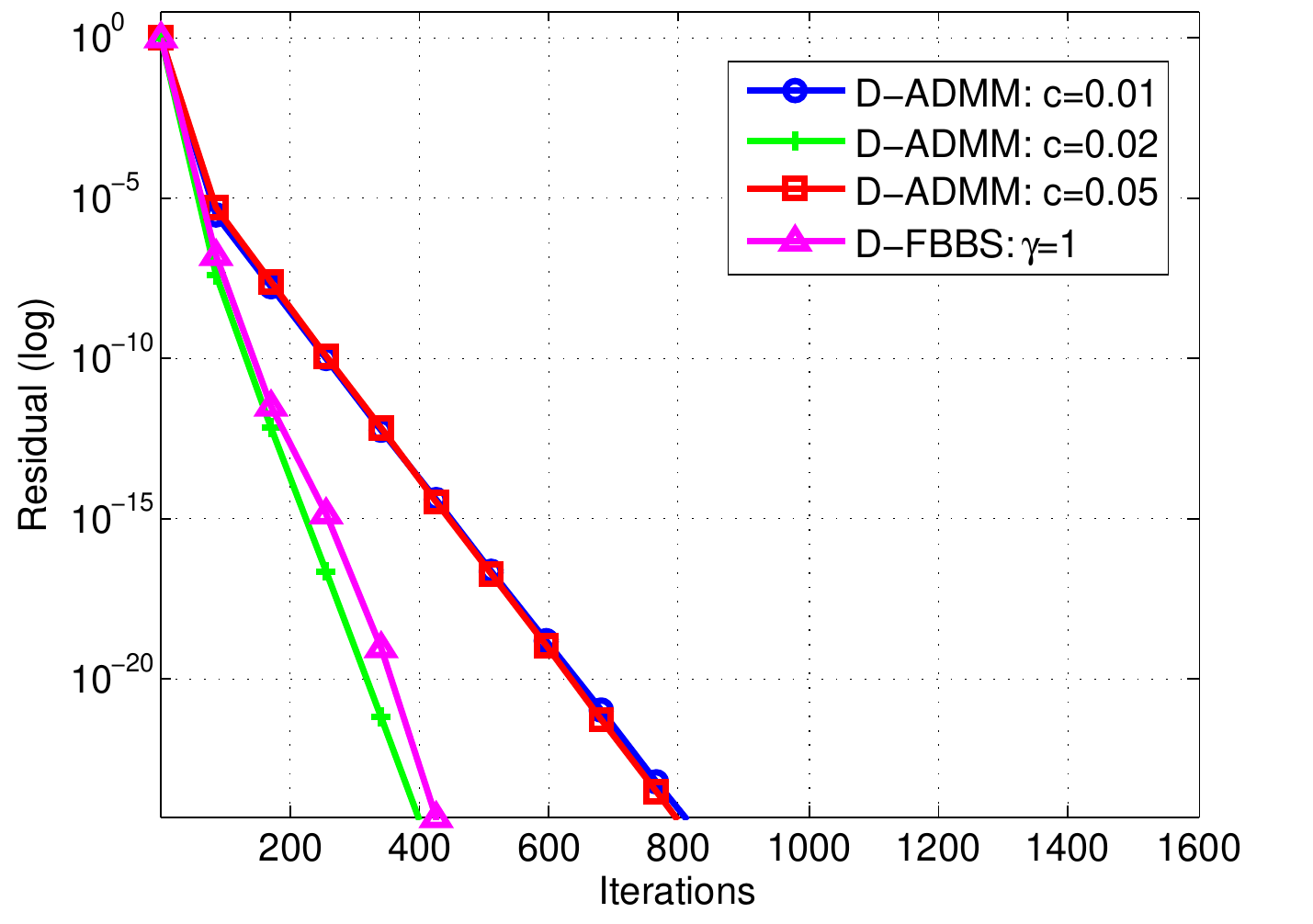}\label{fig:slowcomm}}
 \subfloat[$r=2$]{ \includegraphics[width=0.32\textwidth]{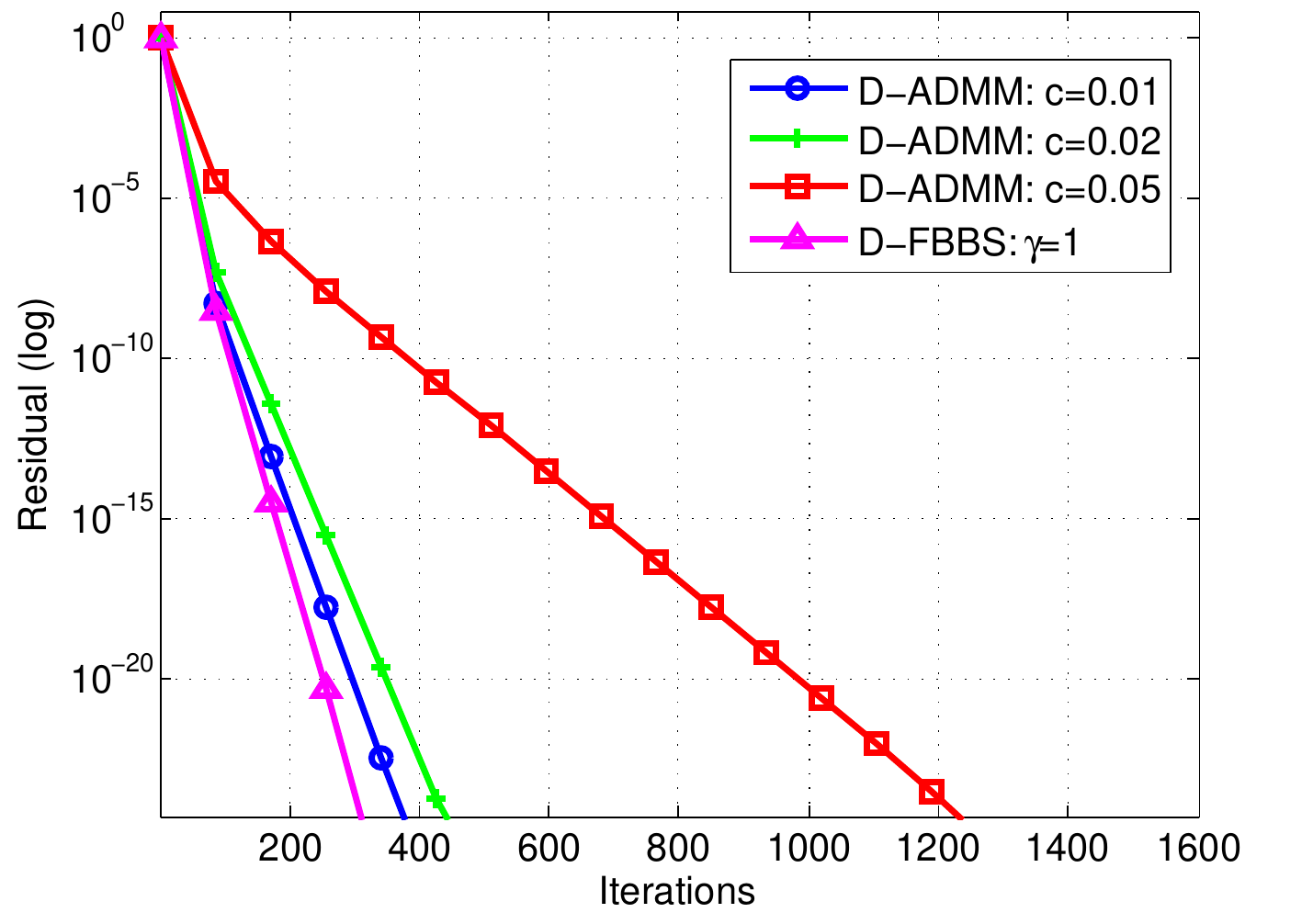}\label{fig:fastcomm}}
\caption{\alert{Performance comparison of the non-gradient-based distributed algorithms: D-ADMM (C-ADMM) and D-FBBS. The optimal parameter values for D-ADMM under the different sensing ranges: $r=1,1.5$ and $2$ are $c=0.05,0.02$ and $0.01$ respectively.}}
\label{fig:algcomparison_DFBBS_DADMM}
\end{figure*}

\begin{figure*}[!htpb]
\centering
 \subfloat[$r=1$]{ \includegraphics[width=0.32\textwidth]{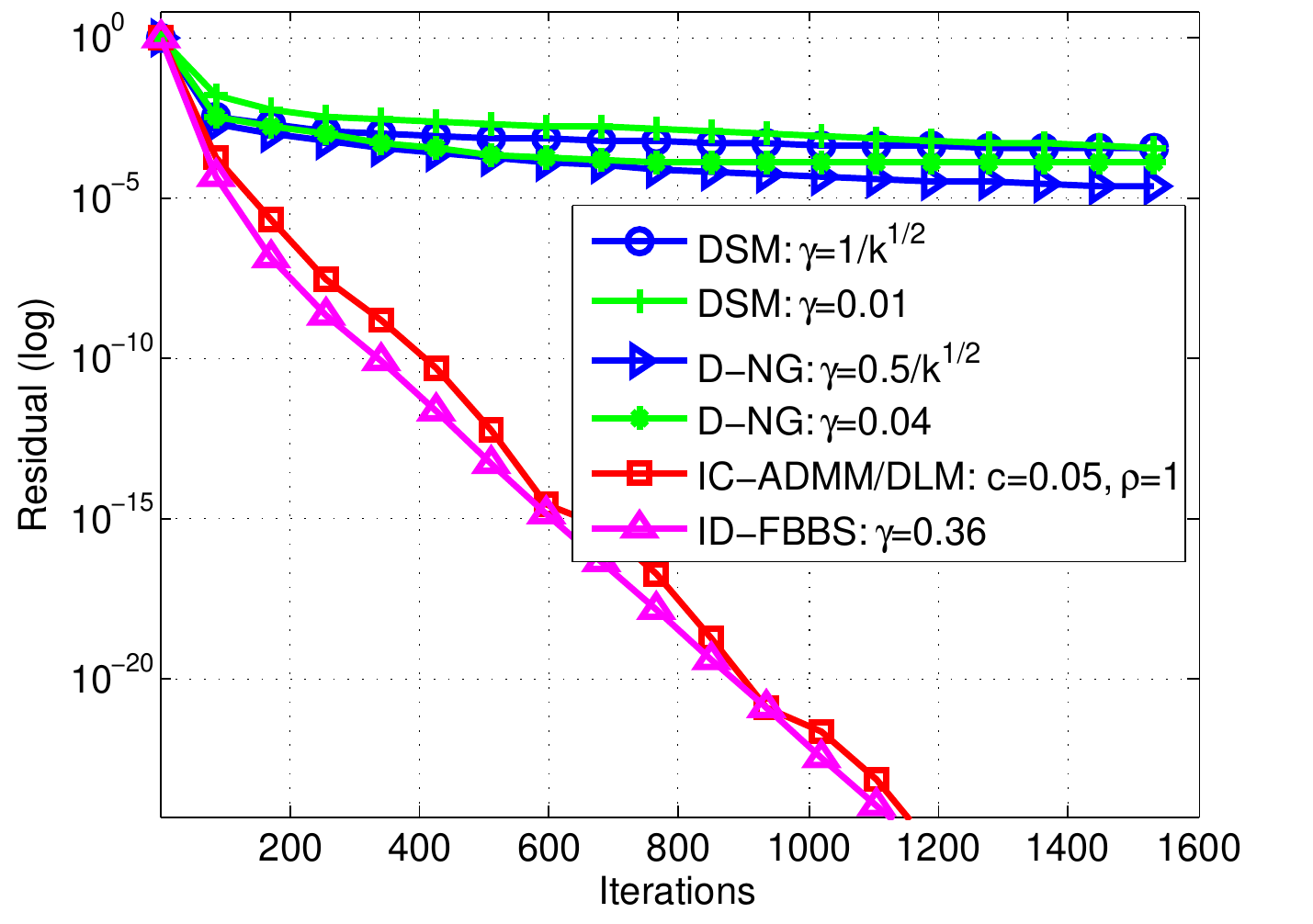}\label{fig:nocomm}}
 \subfloat[$r=1.5$]{ \includegraphics[width=0.32\textwidth]{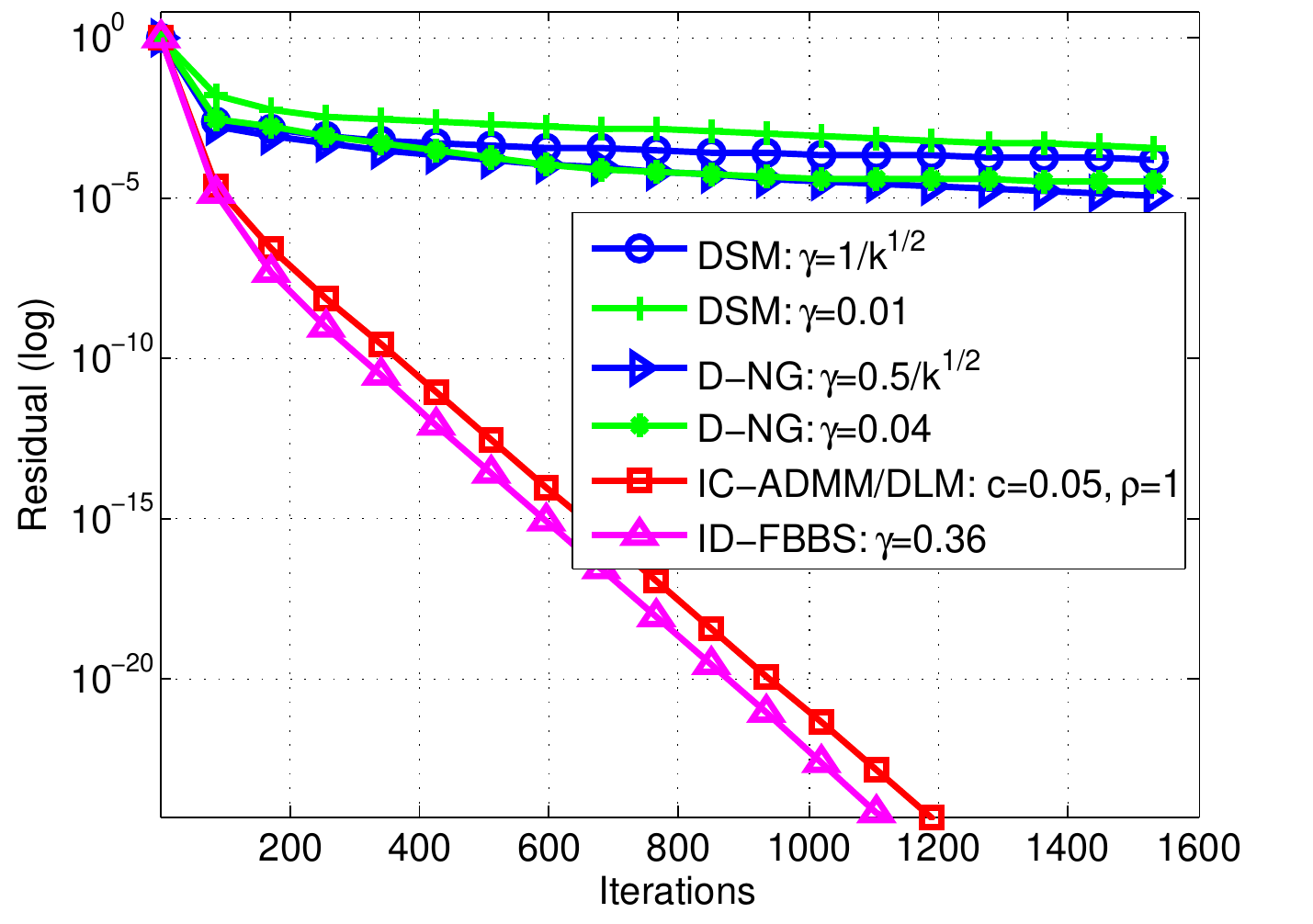}\label{fig:slowcomm}}
 \subfloat[$r=2$]{ \includegraphics[width=0.32\textwidth]{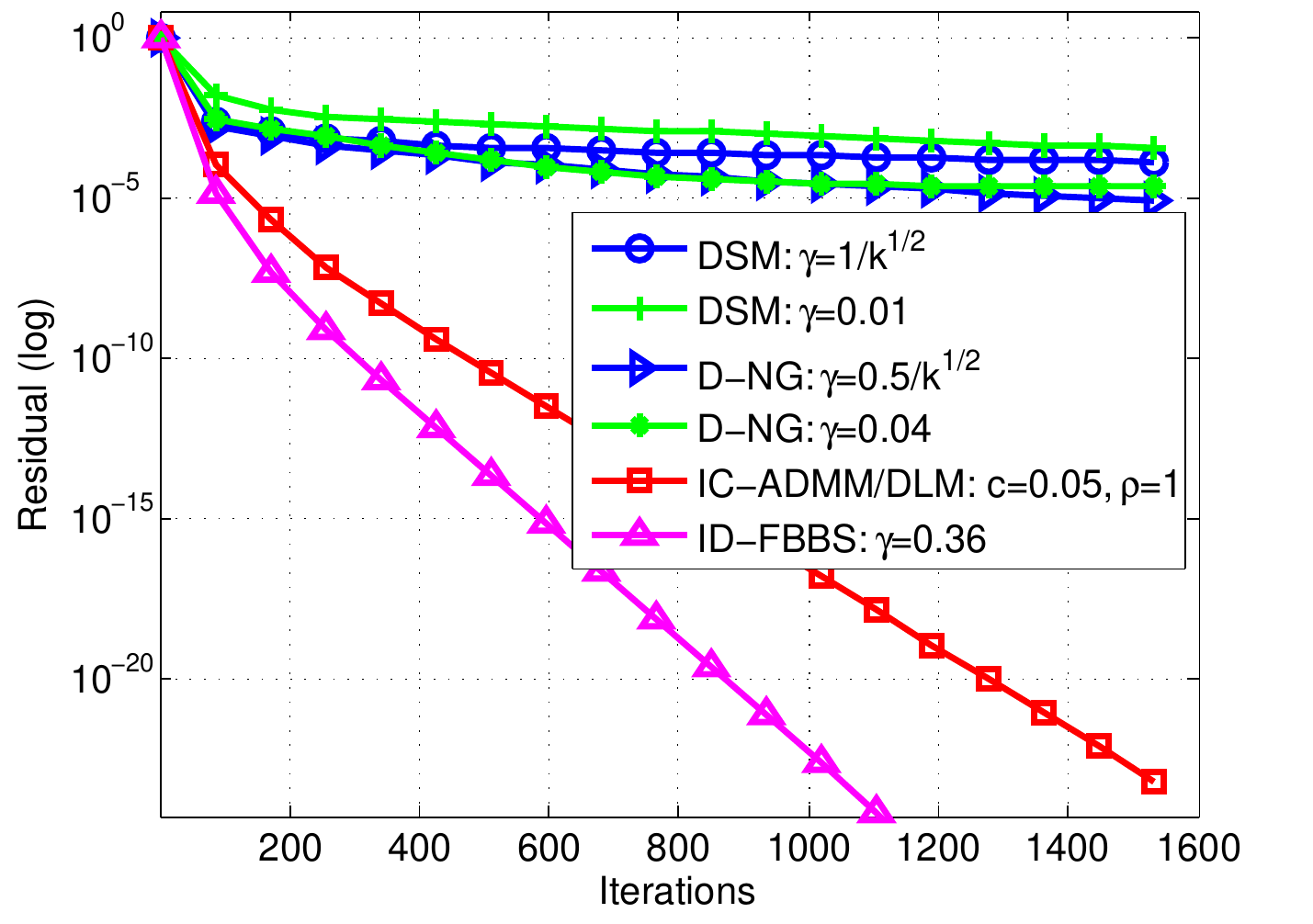}\label{fig:fastcomm}}
\caption{\alert{Performance comparison of the gradient-based distributed algorithms: DSM, D-NG, IC-ADMM/DLM and ID-FBBS. The stepsizes or parameters for each algorithm are hand-optimized for the scenario with a sensing range of $r=1$.}} 
\label{fig:algcomparison_DFBBS_ICADMMorDLM_DSM}
\end{figure*}

{\bf Parameter Setting:} We set $d=4,s=1, m=50$ for both networks. The measurement matrix is generated from a uniform distribution in the unit $\mathcal{R}^{s\times d}$ space and the noise follows a i.i.d. Gaussian process with zero mean and certain variance, i.e., $\mathcal{N}(0,0.1)$.  The regularization parameter is set as $\lambda_{reg}=0.1$ for the stochastic network. For the strong convexity of the cost function is not required for fixed network, we set $\lambda_{reg}=0$ for this case. The weight matrix is designed using the modified\footnote{This is to ensure the positivenss of the weight matrix.} Metropolis-Hastings protocol~\cite{xiao2004fast} as follows:
\begin{equation}\label{eq:metropolis-hastings_weight}
w_{ij}=
\begin{cases}
\frac{1}{2\cdot\max \{d_i,d_j\}},~~~~~~~\text{if}~ (i,j)\in\mathcal{E}\\
1-\sum_{j\in\mathcal{N}_i}w_{ij},~~~\text{if}~ i=j\\
0, ~~~~~~~~~~~~~~~~~~~\text{otherwise,}
\end{cases}
\end{equation}
where $d_i,d_j$ are the degrees of node $i$ and $j$ receptively.
The simulation results for the stochastic network are averaged over 20 runs. \alert{For the sake of fairness, we compare our D-FBBS algorithm with D-ADMM/C-ADMM~\cite{shi2014linear,mateos2010distributed} which is also a non-gradient-based approach while comparing our inexact algorithm (ID-FBBS) with DSM~\cite{Subgradient_Consensus_Nedic}, D-NG (the accelerate version of DSM)~\cite{Jakovetic2014_FastGradMethod} and IC-ADMM/DLM~\cite{chang2015multi,ling2015dlm} which are all gradient-based approaches thus having comparable computational complexity. For stochastic networks, we compare our algorithm D-FBBS with DSM which is also applicable to the same stochastic settings as this work. All comparisons are made in terms of the relative fixed point residual $\text{FPR}=\frac{\norm{x_k-x^\star}^2}{\norm{x_0-x^\star}^2}$ for both networks.} 

{\bf Discussions:} 
\alert{Fig.~\ref{fig:algcomparison_DFBBS_DADMM} plots the relative FPR of D-FBBS and D-ADMM/C-ADMM with respect to the number of iterations for a fixed networks under certain parameter setting. The optimal parameter setting for D-ADMM/C-ADMM under the graph with the different sensing ranges (i.e., $r=1,~1.5$ and $2$) is $c=0.05,~0.02$ and $0.01$ respectively. It follows from the figure that the value of the parameter $c$ that is optimal to one particular graph (e.g., $r=1$) may not be suitable to others (e.g., $r=2$). For instance, $c=0.05$ is optimal for the graph with $r=1$ but the same value of the parameter will result in a dramatic performance degradation for the graph with $r=2$. In contrast, our algorithm has a comparable performance under the same parameter stetting of $\gamma=1$ regardless of the sensing range and the advantage is more significant when the network is more connected. In other words, the proposed algorithm is, to some extent, more robust to the variation of the network, implying its capacity in dealing with varying networks.}

\begin{figure}[!hbpt] 
      \centering
      \includegraphics[scale=0.54]{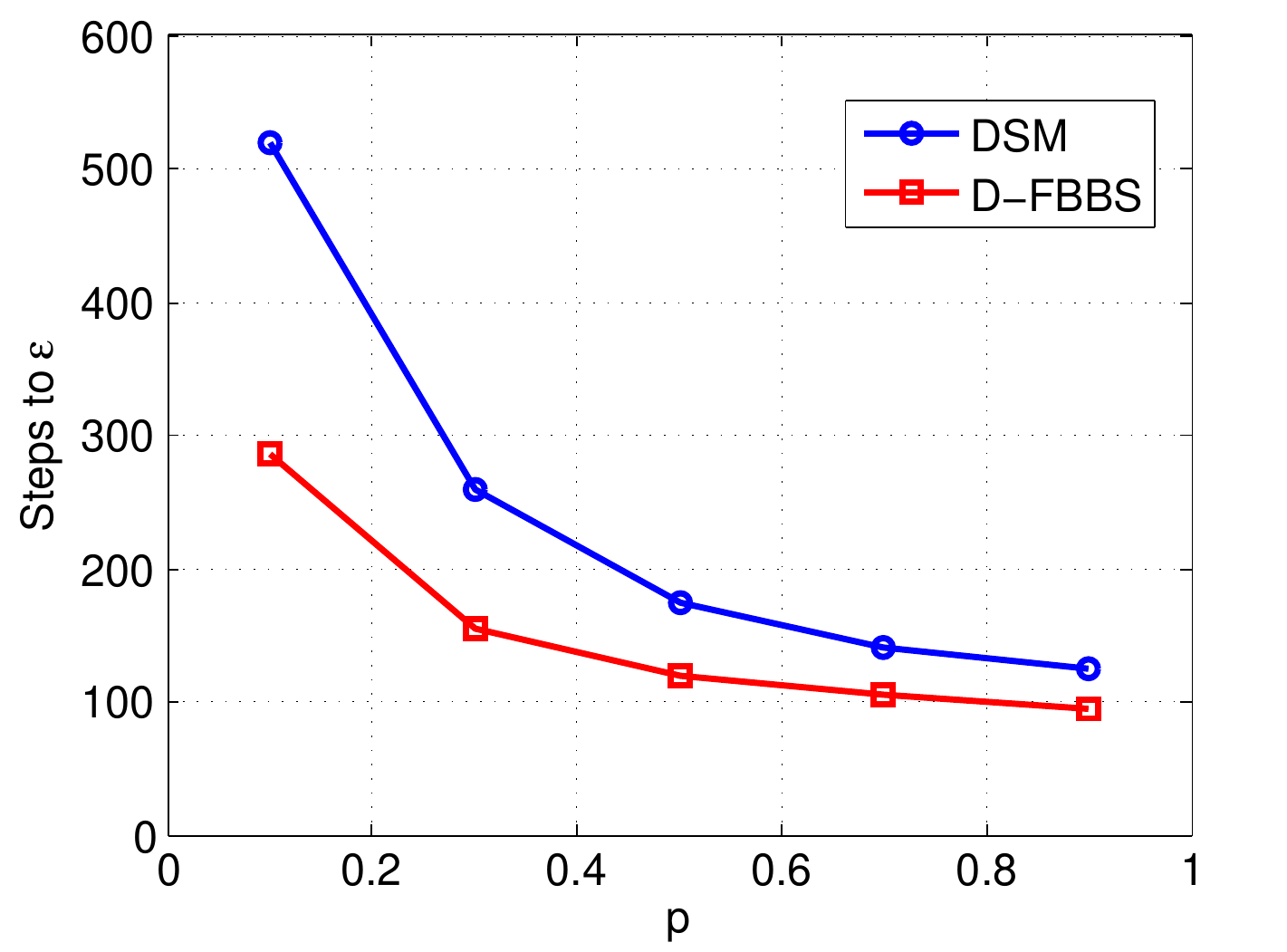}   
      \caption{Plot of the number of iterations required to reach a fixed accuracy $\epsilon=0.001$ for both DSM and D-FBBS. The stepsize $\gamma=2/k$ for DSM is optimized by hand while the stepsize $\gamma=10$ for D-FBBS is calculated based on Theorem~\ref{thm:convergence_mainresult_rndnet}. The results are obtained over an average of 20 runs.}
      \label{fig:splitBregman_vs_DSM_stepsbyplink}
\end{figure}
\begin{figure}[!hbpt] 
      \centering
      \includegraphics[scale=0.54]{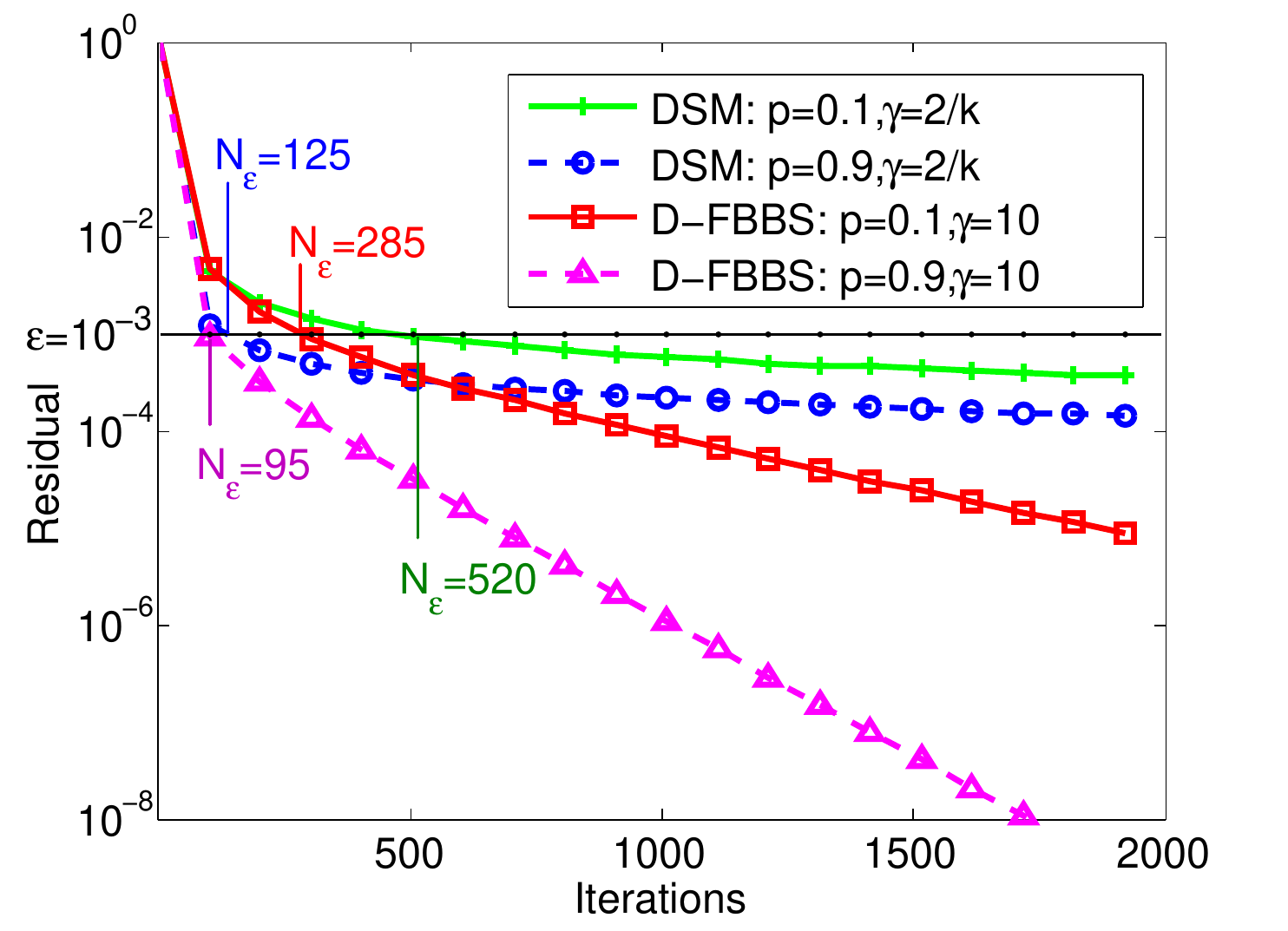}   
      \caption{Plot of the relative FPR versus the number of iterations for both DSM and D-FBBS under two different  probabilities of link failure, i.e., $p=0.1$ (high probability of failure) and $p=0.9$ (low probability of failure).}
      \label{fig:splitBregman_vs_DSM_rndnet}
\end{figure}

\alert{Fig.~\ref{fig:algcomparison_DFBBS_ICADMMorDLM_DSM} plots the relative FPR of the gradient-based algorithms, DSM, D-NG, IC-ADMM/DLM and ID-FBBS,  with respect to the number of iterations for a fixed network under certain stepsizes. It follows from the figure that all algorithms have similar convergence performance at the initial stage. However, DSM and D-NG either with constant stepsize or decaying stepsize get stuck after several iterations, which is not surprising as implied by the theoretical result in~\cite{yuan2013convergence} and the slow convergence rate of using decaying stepsize. In contrast, the proposed ID-FBBS algorithm as well as IC-ADMM/DLM still converge almost linearly to the optimum. Additionally, similar with the above analysis, it can be observed from the figure that the convergence performance of IC-ADMM/DLM, under the parameter setting of $c=0.05,\rho=1$ that are optimized for the graph with $r=1$, is varying (becoming worse) with the sensing range. Instead, our algorithm, under the same constant stepsize $\gamma=0.36$, perform almost the same for all scenarios with different sensing ranges. }

Fig.~\ref{fig:splitBregman_vs_DSM_stepsbyplink} and Fig.~\ref{fig:splitBregman_vs_DSM_rndnet} illustrate the simulation results for stochastic network. Fig.~\ref{fig:splitBregman_vs_DSM_stepsbyplink} shows that the proposed D-FBBS algorithm always needs less iterations to reach the specified accuracy of $\epsilon=0.001$ as compared with DSM. The advantage is more significant when the communication link has low probability of being activated at each iteration, i.e., communication process being more ``asynchronous". In particular, we can observe from Fig.~\ref{fig:splitBregman_vs_DSM_rndnet} that the proposed D-FBBS algorithm require $N_\epsilon=95$ iterations to reach the specified accuracy while DSM needs $N_\epsilon=125$ iterations when the communication link is activated with a high probability of $p=0.9$. In addition, when $p$ becomes small (e.g., $p=0.1$), the difference of iterations required is enlarged, i.e., $N_\epsilon=285$ for D-FBBS and $N_\epsilon=520$ for DSM. Moreover, similar to the case of fixed network, both algorithms have similar performance in the beginning but the proposed D-FBBS algorithm still progress at an almost linear convergence rate afterwards.

\section{Conclusion}
In this paper, we proposed a new distributed algorithm framework based on Bregman method and operator splitting, which allows us to recover most existing well-known distributed algorithms. In particular, based on Bregman iterative regularization and forward-backward splitting technique, we developed a new algorithm (D-FBBS) for distributed optimization problems. We have showed that, for general convex cost functions, the  algorithm are able to achieve a non-ergodic convergence rate of $o(\frac{1}{k})$ for fixed networks while, for stochastic networks, an ergodic convergence rate of $O(\frac{1}{k})$ can be obtained for strongly convex cost functions. In addition, an inexact version of D-FBBS (ID-FBBS) has also been provided for more computational efficiency and shown to also have a non-ergodic convergence rate of $o(1/k)$. Moreover, several simulations were done to verify the advantages of the proposed algorithms compared to the existing ones.

However, there are several interesting issues to be explored in the future work. For instance, since our result for stochastic networks (when $p=1$) is not ``compatible'' with that of fixed networks, it would be of interest to model explicitly the changing of the topology, establishing stronger convergence result that includes fixed networks as a special case. In addition, it would be worthy of investigation to consider the scenario where the computation is asynchronous as well.


%

\appendices
\section{}\label{sec:appendix}
\subsection{Proof of Proposition~\ref{prop:optim_iff_saddlepoint}}
\begin{IEEEproof}
Since $\ones^Ty=0$ and $null(I-W)=span\{\ones\}$, together with (\ref{eq:dops_optimality_conditions_primalfeasible}) we have $\psi(x^\star,y)=f(x^\star)-y^Tx^\star+\frac{1}{2\gamma}\norm{x^\star}^2_{I-W}=f(x^\star)=\psi(x^\star,y^\star)$. Thus, the left-hand side of condition~(\ref{eq:saddlepoint_definition}) is proved. Then, from (\ref{eq:dops_optimality_conditions_primalfeasible}) and (\ref{eq:dops_optimality_conditions_Lagrangeoptimal}), we have $0\in\partial f(x^\star)+\frac{1}{\gamma}(I-W)x^\star-y^\star=\partial \psi_x(x^\star,y^\star)$, which implies that $x^\star$ minimizes $\psi(x,y^\star)$. Thus, we prove the right-hand side of condition~(\ref{eq:saddlepoint_definition}). Conversely, assume that $(x^\star,y^\star)$ is a saddle point such that the condition~(\ref{eq:saddlepoint_definition}) holds. Then, from the left-hand side of condition~(\ref{eq:saddlepoint_definition}) we have 
\begin{equation*}
\begin{aligned}
\sup_{\ones^Ty=0}\psi(x^\star,y)&=\sup_{\ones^Ty=0}f(x^\star)-y^Tx^\star+\frac{1}{2\gamma}\norm{x^\star}^2_{I-W}\\
&=\psi(x^\star,y^\star)<\infty,
\end{aligned}
\end{equation*}
which is only possible when $x^\star\in span\{\ones\}$. Thus we have $x^\star\in null\{I-W\}$ or namely $(I-W)x^\star=0$.
In addtion, from the righthand side of condition~(\ref{eq:saddlepoint_definition}), we know that $x^\star$ minimizes $\psi(x,y^\star)$, implying that $0\in\partial f(x^\star)+\frac{1}{\gamma}(I-W)x^\star-y^\star$. Since we have shown that $(I-W)x^\star=0$, we have $0\in \partial f(x^\star)-y^\star$, i.e., $y^\star\in\partial f(x^\star)$. Since $y^\star\in span^\perp\{\ones\}$ and thus $\ones^Ty^\star=0$, all the optimality conditions~(\ref{eq:dops_optimality_conditions}) hold. Moreover, since the above analysis shows that every saddle point satisfies the optimality conditions~(\ref{eq:dops_optimality_conditions}), it follows from~\cite[Th. 19.1]{bauschke2011convex} that $x^\star$ is a primal solution to the OCP problem and $y^\star$ is a dual solution to the OEP problem respectively.
\end{IEEEproof}
\subsection{Proof of Lemma~\ref{lem:bijective_transform}}
\begin{IEEEproof}
Since $null(P)= span\{\ones\}$, we have $rank(P)=m-1$. In addition, ${y}\in span^\perp\{\ones\}$ implies that $\ones^T{y}=0$ and in turn that $rank([P~{y}])=m-1$. Thus, by basic linear algebra, there exists a solution ${y}'$ such that $P{y}'={y}$. Since ${y}'\in span^\perp\{\ones\}$ implying that $\ones^T{y}'=0$, we have an augmented system of equation $[P^T~\ones]^T{y}'=[{y}~0]^T$. Since $rank([P^T~\ones])=m$, again by basic linear algebra, we conclude that ${y}'$ is unique. The reverse is similar.
\end{IEEEproof}

\subsection{Proof of Proposition~\ref{prop:convergence_IDFBBS}}
\begin{IEEEproof} 
\alert{ Let $\mathcal{S}$ be the optimal solution set of the \emph{primal-dual} problem and $(x^\star,y^\star)\in\mathcal{S}$ be a saddle point. Define $\Delta f_k=\nabla f(x_k)-\nabla f(x_{k-1}),\Delta x_{k}=x_{k}-x_{k-1}$. Then, adding (\ref{alg:dops_fixednet_ineaxct_Uzawa_smooth_forbacksplit_x}) with (\ref{alg:dops_fixednet_ineaxct_Uzawa_smooth_forbacksplit_y}) and rearranging terms yields:
\begin{subequations}\label{alg:dops_fixednet_ineaxct_Uzawa_smooth_forbacksplit_equi}
\begin{align}
\gamma(\Delta f_{k+1}+y_{k+1})-W\Delta x_{k+1}&=\gamma\nabla f(x_{k+1}), \label{alg:dops_fixednet_ineaxct_Uzawa_smooth_forbacksplit_equi_x}\\
(I-W)x_{k+1}+\gamma(y_{k+1}-y_k)&=0.  \label{alg:dops_fixednet_ineaxct_Uzawa_smooth_forbacksplit_equi_y}
\end{align}
\end{subequations} 
Since $\ones^Ty_k=0, \forall k\geq 0$ (cf.~Lem.~\ref{lem:ConservProp}) and $null(I-W)=span\{\ones\}$, it follows from Lemma~\ref{lem:bijective_transform} that there exists a unique $y'_k\in span^\perp\{\ones\}$ such that $y_k=(I-W)y'_k$. Knowing that $\nabla f$ is maximally monotone (cf.~Assum.~\ref{assum:costfuctions_proper_closed_convex}), $\gamma y^\star=\gamma\nabla f(x^\star)$ from the optimality condition~(\ref{eq:dops_optimality_conditions_Lagrangeoptimal}) and using (\ref{alg:dops_fixednet_ineaxct_Uzawa_smooth_forbacksplit_equi_x}),   we have
\begin{equation}\label{eq:fixednet_mainresult_inexact_iterseq_fixedpointresidual_primitive}
\begin{aligned}
&\innprod{\gamma(y_{k+1}-y^\star)-W\Delta x_{k+1}+\gamma\Delta f_{k+1}}{x_{k+1}-x^\star}\\
&=\innprod{\gamma(I-W)(y'_{k+1}-y'^\star)}{x_{k+1}-x^\star}\\
&~~~+\innprod{\gamma\Delta f_{k+1}}{x_{k+1}-x^\star}-\innprod{W\Delta x_{k+1}}{x_{k+1}-x^\star}\geq 0.\\
\end{aligned}
\end{equation}
Let us now consider the first term of the last inequality. Similar with~(\ref{eq:fixednet_mainresult_iterseq_fixedpointresidual_raw}) we can easily derive the following relation:
\begin{equation}
\begin{aligned}
&\innprod{\gamma(I-W)(y'_{k+1}-y'^\star)}{x_{k+1}-x^\star}\\
&~~~~~~~~~=-\gamma^2\innprod{(I-W)(y'_{k+1}-y'_k)}{y'_{k+1}-y'^\star}.
\end{aligned}
\end{equation}
Then, using the similar identity as~(\ref{eq:fixednet_iterseq_identity}) and the identity of Bregman distance (cf.~Eq. (\ref{eq:Bregman_threeidentity})), (\ref{eq:fixednet_mainresult_inexact_iterseq_fixedpointresidual_primitive}) becomes
\begin{equation}\label{eq:fixednet_mainresult_inexact_iterseq_fixedpointresidual_main}
\begin{aligned}
&\gamma^2\norm{y'_{k+1}-y'^\star}^2_{I-W}-\gamma^2\norm{y'_k-y'^\star}^2_{I-W}+\norm{x_{k+1}-x^\star}^2_W\\
&-\norm{x_k-x^\star}^2_W-2\gamma D_f(x^\star,x_{k+1})+2\gamma D_f(x^\star,x_{k})\\
&\leq -\gamma^2\norm{y'_{k+1}-y'_k}^2_{I-W}-\norm{\Delta x_{k+1}}^2_W+2\gamma D_f(x_{k+1},x_k).
\end{aligned}
\end{equation}
Now, recalling that $\tilde{x}_k=(I-\avector)x_k$, $y_{k+1}=(I-W)y'_{k+1}$ and $y_k=(I-W)y'_k$, similar with (\ref{eq:fixednet_mainresult_coordinatetransform_disagreement}) we have  
\begin{equation*}
\begin{aligned}
\gamma^2\norm{y'_{k+1}-y'_k}^2_{I-W}=\norm{\tilde{x}_{k+1}}^2_{I-W}.
\end{aligned}
\end{equation*}
Also, as before, since $\ones^Ty'_k=0, \forall k\geq 0$ (cf.~Lem.~\ref{lem:ConservProp}) and $\rho\left(W-\avector\right)<1$ by Assumption~\ref{assum:fixednet_weightmatrix}, we have $y_{k+1}-y^\star=(I-W)(y'_{k+1}-y'^\star)=[I-(W-\avector)](y'_{k+1}-y'^\star)$ and $$\norm{y'_{k+1}-y'^\star}^2_{I-W}=\norm{y_{k+1}-y^\star}^2_{[I-(W-\avector)]^{-1}}.$$
Let $V_k=\gamma^2\norm{y_{k+1}-y^\star}^2_{[I-(W-\avector)]^{-1}}+\norm{x_k-x^\star}^2_W-2\gamma D_f(x^\star,x_{k})$, which is positive since $\gamma<\frac{\lambda_{\min}(W)}{L_f}$ such that $\norm{x_k-x^\star}^2_W\geq\lambda_{\min}(W)\norm{x_k-x^\star}^2>\gamma L_f\norm{x_k-x^\star}^2\geq2\gamma D_f(x^\star,x_{k})$. Thus, (\ref{eq:fixednet_mainresult_inexact_iterseq_fixedpointresidual_main}) can be rewritten as follows
\begin{equation}\label{eq:fixednet_mainresult_inexact_iterseq_fixedpointresidual_main_equi}
\begin{aligned}
&V_{k}-V_{k+1}\\
&~~~\geq \norm{\tilde{x}_{k+1}}^2_{I-W}+\norm{\Delta x_{k+1}}^2_W-2\gamma D_f(x_{k+1},x_k).
\end{aligned}
\end{equation}
Summing (\ref{eq:fixednet_mainresult_inexact_iterseq_fixedpointresidual_main_equi}) over $k$ from $0$ to $n-1$ and knowing that all eigenvalues $\lambda (W)\in (0,1]$ by Assum.~\ref{assum:fixednet_weightmatrix} (cf. Rem.~\ref{rem:weightmatrix_eigen_speed}) leads to
\begin{equation}
\begin{aligned}
\infty&>V_0-V_n\\
&\geq\sum_{k=0}^{n-1}\left(\norm{\tilde{x}_{k+1}}^2_{I-W}+\norm{\Delta x_{k+1}}^2_W-2\gamma D_f(x_{k+1},x_k)\right)\\
&\geq\sum_{k=0}^{n-1}\left(\norm{\tilde{x}_{k+1}}^2_{I-W}+\norm{\Delta x_{k+1}}^2_{\lambda_{\min}(W)-\gamma L_f}\right)\\
&\geq\bracket{\lambda_{\min}(W)-\gamma L_f}\sum_{k=0}^{n-1}\bracket{\norm{\tilde{x}_{k+1}}^2_{I-W}+\norm{\Delta x_{k+1}}^2_W}.
\end{aligned}
\end{equation}
From Lemma~\ref{lem:fixednet_fixedpointresidual_non_increasing}, we know that $u_k=\norm{\tilde{x}_{k}}^2_{I-W}+\norm{\Delta x_{k}}^2_W$ is monotonically non-increasing. Thus, we have
\begin{equation}
\begin{aligned}
&nu_n\leq\sum_{k=0}^{n-1}u_{k+1}\leq \frac{V_0}{\lambda_{\min}(W)-\gamma L_f},
\end{aligned}
\end{equation}
yielding $u_n\leq o(\frac{1}{n})$.
Since $u_n\geq 0$ by Assumption~\ref{assum:fixednet_weightmatrix}, we claim that $\lim_{n\rightarrow\infty}u_n=0$. From (\ref{eq:fixednet_mainresult_inexact_iterseq_fixedpointresidual_main_equi}) we know that $V_k$ is bounded and so is $x_k$ and $y_k$. Using the similar argument as in the proof of Theorem~\ref{thm:convergence_mainresult_fixednet}, it follows that $\seq{(x_k,y_k)}$ has a unique cluster point in the optimal set $\mathcal{S}$. Thus, we conclude that the sequence $\seq{(x_k,y_k)}$ will converge to a saddle point $(x^\star,y^\star)$  with a non-ergodic rate of $o(\frac{1}{k})$.
}
\end{IEEEproof}




\ifCLASSOPTIONcaptionsoff
  \newpage
\fi



\bibliographystyle{IEEEtran}
\bibliography{IEEEfull,reference,ref_splitbregman,ref_convexoptim}
\end{document}